\documentclass[10pt]{article}
\usepackage{amsmath, amsthm, amssymb, enumerate}
\usepackage{color}
\usepackage{datetime}
\usepackage{fancyhdr}

  \chardef\forshowkeys=0
  \chardef\showllabel=0
  \chardef\refcheck=0
  \chardef\sketches=0
  \chardef\showcolors=0
  
\usepackage{enumitem}
\usepackage{datetime}
\usepackage{fancyhdr}
\usepackage{comment}
\allowdisplaybreaks
\ifnum\forshowkeys=1
  
  \usepackage[notref,notcite,color]{showkeys}
\fi
\ifnum\showllabel=1
 \def\llabel#1{\marginnote{\color{lightgray}\rm\small(#1)}[-0.0cm]\notag}
\else
 \def\llabel#1{\notag}
\fi
\usepackage[margin=1in]{geometry}
\usepackage{amsmath, amsthm, amssymb}
\usepackage{times}
\usepackage{graphicx}
\usepackage[usenames,dvipsnames,svgnames,table]{xcolor}
\usepackage{marginnote}
\usepackage[unicode,breaklinks=true,colorlinks=true,linkcolor=blue,urlcolor=blue,citecolor=blue]{hyperref}
\usepackage[most]{tcolorbox}
\usepackage{tikz}

\ifnum\refcheck=1
  \usepackage{refcheck}
\fi

\chardef\coloryes=0 
\chardef\isitdraft=0

\ifnum\isitdraft=1 
\usepackage[color,notcite]{showkeys}
   \def\eqref#1{({\ref{#1}})}                

\definecolor{refkey}{rgb}{.3,0.3,0.3}

\else
  \definecolor{refkey}{rgb}{.8,0.8,0.1}
  \definecolor{labelkey}{rgb}{.9,0.6,0.1}

  \def\startnewsection#1#2{\section{#1}\label{#2}\setcounter{equation}{0}}   
\fi

\begin{document}
\def\BM{\text{BM}}
\def\Deltah{\Delta_{\text h}}  
\def\chere{\tdot\nts{CONTINUE HERE}\tdot}
\def\bdot{{\color{blue} {\hskip-.0truecm\rule[-1mm]{4mm}{4mm}\hskip.2truecm}}\hskip-.3truecm}
\newcommand\refer{\eqref}
\def\nablaa{\nabla_{a}}
\def\lec{\lesssim}
\def\vv{\tilde v}

\ifnum\showcolors=1
  \def\colr{\color{red}}
  \def\colrr{\color{black}}
  \def\colb{\color{black}}
  \def\coly{\color{lightgray}}
  \definecolor{colorgggg}{rgb}{0.1,0.5,0.3}
  \definecolor{colorllll}{rgb}{0.0,0.7,0.0}
  \definecolor{colorhhhh}{rgb}{0.3,0.75,0.4}
  \definecolor{colorpppp}{rgb}{0.7,0.0,0.2}
  \definecolor{coloroooo}{rgb}{0.45,0.0,0.0}
  \definecolor{colorqqqq}{rgb}{0.1,0.7,0}
  \def\colg{\color{colorgggg}}
  \def\collg{\color{colorllll}}
  \def\cole{\color{coloroooo}}
  \def\coleo{\color{colorpppp}}
  \def\colu{\color{blue}}
  \def\colc{\color{colorhhhh}}
  \def\colW{\colb}   
  \definecolor{coloraaaa}{rgb}{0.6,0.6,0.6}
  \def\colw{\color{coloraaaa}}
\else
  \def\colr{\color{black}}
  \def\colrr{\color{black}}
  \def\colb{\color{black}}
  \def\coly{\color{black}}
  \def\colg{\color{black}}
  \def\collg{\color{black}}
  \def\cole{\color{black}}
  \def\coleo{\color{black}}
  \def\colu{\color{black}}
  \def\colc{\color{black}}
  \def\colW{\color{black}}
  \def\colw{\color{black}}
\fi

\def\bnew{\colr {\bf NEW:}}
\def\enew{\colb {}}
\def\bold{\colu {\bf OLD:}}
\def\eold{\colb {}}
\def\largesign#1{\colu {\text{\huge #1}}}
\def\lot{\text{l.o.t.}}

\def\inon#1{\hbox{\ \ \ \ \ \ \ }\hbox{#1}}                
\def\onon#1{\inon{on~$#1$}}
\def\inin#1{\inon{in~$#1$}}

\def\Dg{{D'g}}
\def\ua{u^{\alpha}}
\def\into{\int_{\Omega}}
\def\ques{{\cor \underline{??????}\cob}}
\def\nto#1{{\coC \footnote{\em \coC #1}}}
\def\fractext#1#2{{#1}/{#2}}
\def\fracsm#1#2{{\textstyle{\frac{#1}{#2}}}}   
\def\baru{U}
\def\nnonumber{}
\def\palpha{p_{\alpha}}
\def\valpha{v_{\alpha}}
\def\qalpha{q_{\alpha}}
\def\walpha{w_{\alpha}}
\def\falpha{f_{\alpha}}
\def\dalpha{d_{\alpha}}
\def\galpha{g_{\alpha}}
\def\halpha{h_{\alpha}}
\def\psialpha{\psi_{\alpha}}
\def\psibeta{\psi_{\beta}}
\def\betaalpha{\beta_{\alpha}}
\def\gammaalpha{\gamma_{\alpha}}
\def\Talpha{T}
\def\TTalpha{T_{\alpha}}
\def\TTalphak{T_{\alpha,k}}
\def\falphak{f^{k}_{\alpha}}
\def\R{\mathbb R}

\newcommand {\Dn}[1]{\frac{\partial #1  }{\partial N}}
\def\mm{m}

\def\cor{{}}
\def\cog{{}}
\def\cob{{}}
\def\coe{{}}
\def\coA{{}}
\def\coB{{}}
\def\coC{{}}
\def\coD{{}}
\def\coE{{}}
\def\coF{{}}

\ifnum\coloryes=1

  \definecolor{coloraaaa}{rgb}{0.1,0.2,0.8}
  \definecolor{colorbbbb}{rgb}{0.1,0.7,0.1}
  \definecolor{colorcccc}{rgb}{0.8,0.3,0.9}
  \definecolor{colordddd}{rgb}{0.0,.5,0.0}
  \definecolor{coloreeee}{rgb}{0.8,0.3,0.9}
  \definecolor{colorffff}{rgb}{0.8,0.3,0.9}
  \definecolor{colorgggg}{rgb}{0.5,0.0,0.4}
  \definecolor{colorhhhh}{rgb}{0.6,0.6,0.6}

 \def\cog{\color{colordddd}}
 \def\coy{\color{colorhhhh}}
 \def\cogray{\color{colorhhhh}}
 \def\cob{\color{black}}

 \def\cor{\color{red}}
 \def\coe{\color{colorgggg}}

 \def\coA{\color{coloraaaa}}
 \def\coB{\color{colorbbbb}}
 \def\coC{\color{colorcccc}}
 \def\coD{\color{colordddd}}
 \def\coE{\color{coloreeee}}
 \def\coF{\color{colorffff}}
 \def\coG{\color{colorgggg}}

\fi
\ifnum\isitdraft=1
   \chardef\coloryes=1
   \baselineskip=17pt
   \def\blackdot{{\color{red}{\hskip-.0truecm\rule[-1mm]{4mm}{4mm}\hskip.2truecm}}\hskip-.3truecm}
   \def\purpledot{{\coA{\rule[0mm]{4mm}{4mm}}\cob}}
   \def\pdot{\purpledot}
  \definecolor{labelkey}{rgb}{.5,0.1,0.1}
\else
   \baselineskip=15pt
   \newcommand{\bcirc}[2][green,fill=blue]{\tikz[baseline=-0.5ex]\draw[#1,radius=#2] (0,0) circle ;}
   \def\bnew{{\colu {\hskip-.0truecm\rule[-1mm]{4mm}{4mm}\hskip.2truecm}}\hskip-.3truecm}
   \def\blackdot{\color{red}{\rule[-0mm]{4mm}{4mm}}}
   \def\purpledot{{\rule[-3mm]{8mm}{8mm}}}
   \def\pdot{}
\fi

\def\tdot{{\color{green}{\hskip-.0truecm\rule[-1mm]{4mm}{4mm}\hskip.2truecm}}\hskip-.3truecm}
\def\nts#1{{\hbox{\bf ~#1~}}} 
\def\nts#1{{\colu\hbox{\bf ~#1~}}} 
\def\ntsf#1{\footnote{\hbox{\bf ~#1~}}} 
\def\ntsf#1{\footnote{\cor\hbox{\bf ~#1~}}} 
\def\bigline#1{~\\\hskip2truecm~~~~{#1}{#1}{#1}{#1}{#1}{#1}{#1}{#1}{#1}{#1}{#1}{#1}{#1}{#1}{#1}{#1}{#1}{#1}{#1}{#1}{#1}\\}
\def\biglineb{\bigline{$\downarrow\,$ $\downarrow\,$}}
\def\biglinem{\bigline{---}}
\def\biglinee{\bigline{$\uparrow\,$ $\uparrow\,$}}

\def\mbar{{\overline M}}
\def\tilde{\widetilde}
\newtheorem{Theorem}{Theorem}[section]
\newtheorem{Corollary}[Theorem]{Corollary}
\newtheorem{Proposition}[Theorem]{Proposition}
\newtheorem{Lemma}[Theorem]{Lemma}
\newtheorem{Remark}[Theorem]{Remark}
\newtheorem{definition}{Definition}[section]
\def\theequation{\thesection.\arabic{equation}}
\def\endproof{\hfill$\Box$\\}
\def\square{\hfill$\Box$\\}
\def\comma{ {\rm ,\qquad{}} }            
\def\commaone{ {\rm ,\qquad{}} }         
\def\dist{\mathop{\rm dist}\nolimits}    
\def\sgn{\mathop{\rm sgn\,}\nolimits}    
\def\Tr{\mathop{\rm Tr}\nolimits}    
\def\div{\mathop{\rm div}\nolimits}  
\def\TT{R}
\def\supp{\mathop{\rm supp}\nolimits}
\def\divtwo{\mathop{{\rm div}_2\,}\nolimits}
\def\curl{\mathop{\rm curl}\nolimits}    
\def\dbar{\overline\partial}
\def\l{\langle}
\def\r{\rangle}
\def\plusdelta{+\delta}
\def\pd{+\delta}
\def\re{\mathop{\rm {\mathbb R}e}\nolimits}    
\def\indeq{\qquad{}\!\!\!\!}            
\def\period{.}                          
\def\semicolon{\,;}                     
\def\andand{\text{~~~~and~~~~}}
\newcommand{\cD}{\mathcal{D}}

\title{On the interaction between compressible inviscid flow\\and an elastic plate}
\author{ 
Igor Kukavica, \v S\'arka Ne\v casov\'a,
and Amjad Tuffaha} \maketitle
\date{}
\bigskip

\bigskip
\indent Department of Mathematics\\
\indent University of Southern California\\
\indent Los Angeles, USA\\
\indent e-mail: kukavica\char'100usc.edu

\bigskip
\indent Institute of Mathematics\\
\indent Czech Academy of Sciences\\
\indent Prague, Czech Republic\\
\indent e-mail: matus\char'100math.cas.cz

\bigskip
\indent Department of Mathematics\\
\indent American University of Sharjah\\
\indent Sharjah, UAE\\
\indent e-mail: atufaha\char'100aus.edu

\bigskip
\begin{abstract}

We address a free boundary model for the compressible Euler equations where the free boundary, which is elastic, evolves according to a weakly damped fourth order hyperbolic equation forced by the fluid pressure. This system captures the interaction of an inviscid fluid with an elastic plate.  We establish a~priori estimates on local-in-time solutions in low regularity Sobolev spaces, namely with velocity and density initial data $v_{0}, R_{0}$ in~$H^{3}$.  The main new device is a variable coefficients space tangential-time differential operator of order~1 with non-homogeneous boundary conditions, which captures the hyperbolic nature of the compressible Euler equations as well as the coupling with the structural dynamics.

\colb
\end{abstract}

\startnewsection{Introduction}{sec1}

In this paper, we address the compressible Euler equations on a domain with a free moving boundary, which is elastic and moves according to a fourth order linear equation forced by the fluid pressure. The interface is allowed to move in the transversal direction, and thus the interaction between the fluid and the structure is captured by the kinematic condition at the interface. We assume that the fluid is barotropic and that it has a non-vanishing density (a liquid), and we study the problem in three space dimensions.  This paper is the first mathematical treatment of this model. A significant challenge in the study of this inviscid model is the dearth of mathematical literature on the variable coefficients Euler equations with non-homogeneous boundary conditions. The analogous models in case of compressible and incompressible Navier-Stokes are very well-studied in the literature and often rely on maximal regularity theorems with a non-homogeneous divergence and boundary data to study the system. Our system requires a new approach to the existence results for a coupled hyperbolic variable coefficients system.
 
The initial study of compressible Euler equation goes back to Agemi and Beir\~ao Da~Veiga \cite{A, B} and later to Schochet~\cite{Sc}. The first treatment addressing the compressible Euler equations on a time dependent domain was due to P.~Secchi~\cite{S}, who considered the case when a motion of the domain was given.  Well-posedness of the free boundary compressible Euler equations for a liquid with vanishing pressure at the boundary was first established by Lindblad~\cite{Li}. The case of gas where the density could vanish requires a different approach, and well-posedness results were obtained in the 1D case by Jang and Masmoudi~\cite{JM1}, in 3D by Coutand, Lindblad, and Shkoller \cite{CLS} and Jang and Masmoudi~\cite{JM2} with no surface tension, in~\cite{CHS} with surface tension, and more recently in~\cite{Lu}.  More recent treatments of the free boundary compressible Euler equations have also considered the incompressible limit of the compressible Euler equations \cite{LL}, and the minimal regularity assumptions on the initial data for local-in-time well-posedness \cite{DK1,DK2}; see also \cite{CS1, CS2} for the compressible vortex sheet problem; see also \cite{L} for a weighted estimates on the compressible Euler equation in the case of a liquid.

On the other hand, fluid-structure interaction models involving the Navier-Stokes equations in both the incompressible and compressible cases have been extensively treated in various settings with both weak and local-in-time strong solutions with and without structural damping. Early works on weak solutions of the incompressible free boundary model with very strong damping were considered in~\cite{CDEG,CCS}. Other works on the incompressible model with the Koiter shell (both linear and nonlinear) modeling blood flow in the arteries were studied in~\cite{MC1, MC2, MC3,MS}.  Global in time solutions for small data were also obtained for strongly damped structures. More recent works on strong solutions have considered the Gevrey regularity of the solutions to the incompressible viscous model without structural damping~\cite{BT}.

 Local-in-time strong solutions for the viscous compressible model were studied in \cite{MRT, M, MT} while weak solutions were considered by \cite{Bo, BS1, T,TW}. Other more intricate models involving heat exchange have also been considered \cite{MM+}.  However, the inviscid problem requires new tools and a different methodology. In fact, the coupling of the Euler equations involving variable coefficient setting with the plate equation, presents a challenge that requires new tools and a different methodology.  It is notable that the only available well-posedness results for compressible inviscid fluid-structure interaction in the literature treat fixed domains and involve either linearized systems \cite{AGW} or linear potential flow equations (irrotational flow) and nonlinear structural equations~\cite{W, LW}.  For the incompressible inviscid model, we refer the reader to the recent works by two of the authors~\cite{KT1, KT2}.

In this paper, we allow the flow to be rotational and consider a periodic setting in which a channel domain has a fixed bottom and moving top.  The main feature in our results is that the regularity guaranteed by the a~priori estimates provides low regularity results in terms of the norms of initial data. We choose to delegate the construction of solutions to a future work, and focus here on a~priori estimates.
One of the main challenges in the study of the this system concerns the trace of the fluid pressure lacking sufficient regularity needed to study the system using a direct fixed point scheme. Moreover, analysis of the density produces a variable coefficient hyperbolic system, with Neumann type data which requires special trace regularity theorems~\cite{Mi} and their sharp improvements \cite{LT}, adapted to the variable coefficients case, which is a serious issue. It is notable that a direct fixed point scheme for our problem is not possible due to the insufficient regularity of the pressure.

To prove a~priori estimates for the system, we rely on a combination of energy estimates for the variable coefficients Euler-plate system and estimates for a hyperbolic system with variable coefficients and Neumann type conditions satisfied by the density. The derived a~priori estimates are local-in-time and are sharp in the sense that they depend on the critical Sobolev norm of the initial data for the compressible Euler equation. Unlike the treatment in~\cite{S}, where the motion of the domain is prescribed with enough regularity to control the aforementioned traces, the regularity is constrained by the coupling between the plate dynamics forced by the pressure and the density equation, and hence the standard treatment of the system as a symmetrizable hyperbolic system \cite{S,CS1,CS2} does not apply.
We use the ALE transformation to fix the domain. The use of the ALE transformation presents advantages for the treatment of interface-fluid interaction problems \cite{MC1}, but compared to the Lagrangian formulation \cite{L, CLS}, presents a serious challenge in the analysis of the density equation.

The main novelty in this paper involves the treatment of a variable coefficient hyperbolic system with non-homogeneous coupled boundary conditions of the form $Q^{2} g+ \div( f g)=F$ satisfied by the log of the density function $g=\log R$. The operator $Q$ is a first order time-space variable coefficient operator with non-homogeneous Neumann type boundary conditions fully coupled to the interface equation.  The Neumann boundary condition is then computed from the kinematic boundary condition, and is thus expressed in terms of the fluid and plate velocities. A certain degree of structural damping in the structure is then necessary to propagate sufficient regularity through the interface velocity into the Neumann boundary condition of our hyperbolic system satisfied by~$g$. For simplicity, we include a Laplace for the structural damping term, but the estimates can be closed also with a slightly weaker degree of damping.
 
In our approach, problematic boundary terms arise in the energy estimates performed on the solution to the hyperbolic system, and trace regularity results via Fourier analysis developed for second order constant coefficient hyperbolic systems with Neumann boundary conditions are not satisfactory for the derivation of uniform estimates, let alone the construction of approximate solutions.  We then exploit the tangential structure of the operator $Q$ when restricted to the boundary in order to control trace terms, and most importantly isolate the problematic boundary term involving product of Dirichlet Neumann traces of the highest order derivatives. Such problematic term, which would normally be given sense through sharp trace regularity results \cite{Mi,LT} in the case of constant coefficients hyperbolic equations, is alternatively treated through the observation that the term is comparable to the rough pressure boundary term arising in a special energy inequality satisfied by the interface displacement function.

The proof proceeds, in several stages. First, we provide a mix of time and tangential estimates on the system consisting of momentum, continuity and interface equations, to control the time derivatives of highest order. We also derive a crucial energy inequality satisfied by the interface structure displacement which allows the control of the problematic boundary terms.  In the second stage, we derive a hyperbolic system satisfied by the log of the density function with Neumann type conditions. We then proceed to derive special energy estimates for the hyperbolic system, for which the tangential time structure of the second order wave operator when restricted to the boundary becomes essential. The next stage of the paper introduces special elliptic estimates for recovering the full Sobolev regularity of the density variables from the time-space regularity established through the estimates on~$Qg$. In the final stage, we also derive vorticity estimates for the compressible Euler equations as well as velocity divergence estimates in order to control the full Sobolev norm of the velocity variable via a div-curl type lemma.

\startnewsection{Model}{sec2}

We consider a system of partial differential equations modeling the interaction between compressible fluid flow and a plate equation. The model consists of the compressible 
Euler equations defined in an evolving domain $\Omega(t)$ and satisfied by a fluid velocity $u$, the fluid pressure $p$ and the fluid density $\rho$
  \begin{align}
  \begin{split}
   u_{t}  + (u \cdot \nabla)u+\frac{1}{\rho} \nabla p &= 0 \inin{\Omega(t) \times [0,T]} 
   \\
   \rho_{t}+ \div(\rho u) &=0 \inin{\Omega(t) \times[0,T]}
   .
  \end{split}
   \llabel{EQ01}
  \end{align}
The boundary of the domain $\Omega(t)$ consists of a rigid part $\Gamma_{0}$ and a moving part $\Gamma_{1}(t)$ signifying the plate.
The transversal scalar displacement variable $w$ of the plate is modeled by a linear fourth order equation 
  \begin{align}
   w_{tt}  + \Deltah^{2} w - \Deltah w_{t} &= p  ~~&\mbox{on}~~ \Gamma_{1}(0) \times [0,T],
   \llabel{EQ02}
  \end{align}
set on the reference configuration $\Gamma_{1}(0)$,
where $p$ is defined on the dynamic interface~$\Gamma_{1}(t)$. The boundary conditions on the plate could be of clamped or other type, but shall be specified below in the context of geometry of the domain under consideration.
The interaction between two domains is captured by the kinematic boundary condition
  \begin{align}
   w_{t}   = u \cdot \nu(x,t)   \inon{$\Gamma_{1}(t) \times [0,T]$},
   \llabel{EQ03}
  \end{align}
where $\nu$ denotes the dynamic outward normal  to the boundary~$\Gamma_{1}(t)$.
On the rigid part of the boundary, we impose
  \begin{align}
   u \cdot N =0 \inon{on $\Gamma_{0} \times [0,T]$},
   \llabel{EQ04}
  \end{align}
where $N$ denotes the outward unit normal.

The pressure is assumed to be a function of the density, typically via a power law. We consider the liquid case and assume a state equation of the form
  \begin{align}
   p= p(\rho)
   ,
   \label{EQ05}
  \end{align}
where $p$ is smooth. We assume that $p'$ is bounded from below in the range of $\rho$
we are considering.
More precisely, with the condition \eqref{EQ162}, we require
  \begin{equation}
    c_1 \leq   p'(R) \leq c_2
    \onon{[m_0/2,2M_0]}    
    .
   \label{EQ102a}
  \end{equation}
For simplicity, we consider the domain to be a periodic channel
  \begin{align}
   \Omega= \mathbb{T}^{2} \times [0,1]
   \llabel{EQ06}
  \end{align}  
with the plate at the top boundary 
  \begin{align}
    \Gamma_1=\Gamma_{1}(0)= \mathbb{T}^{2} \times \{ 1 \}
    ,
   \llabel{EQ07}
  \end{align}  
while the bottom boundary
  \begin{align}
   \Gamma_0=\Gamma_{0}(0)= \mathbb{T}^{2} \times \{0  \}
   \label{EQ08}
  \end{align}  
represents the rigid portion.
We impose the periodic boundary conditions in the $x_{1}$ and $x_{2}$ directions.

We next introduce the ALE change of variable, which fixes the domain.  
We introduce $\psi$  as the harmonic extension of $1+w$, i.e., the solution to the problem
  \begin{align}
  \begin{split}
   \Delta  \psi &= 0     \inin{\Omega \times [0,T]} \\
   \psi &= 1+ w          \onon{\Gamma_{1} \times [0,T]} \\
   \psi &= 0             \onon{\Gamma_{0} \times [0,T]}. 
  \end{split}
   \label{EQ09}
  \end{align}
Note that since we choose an initial condition $w_{0}=0$, we have $\psi(x,0)=x_{3}$.
More general initial condition for $w_{0}$ can easily be addressed
with only notational changes in the statement and the proofs.

The change of variable which fixes the domain is given by the map $\eta\colon \Omega \to \Omega(t)$
defined by 
  \begin{align} 
   \eta( x_{1}, x_{2}, x_{3},t) = (x_{1}, x_{2},  \psi(x_{1}, x_{2}, x_{3},t),t) 
   .
   \label{EQ10}
  \end{align}
We then introduce the ALE variables $v= u \circ \eta$, $q = p \circ
\eta$, and $R=\rho \circ \eta$, i.e.,
  \begin{align}
    \begin{split}
   v(x_{1}, x_{2}, x_{3},t)= u(x_{1}, x_{2}, \psi( x_{1}, x_{2}, x_{3},t),t)\\
   q(x_{1}, x_{2}, x_{3},t)= p(x_{1}, x_{2}, \psi( x_{1}, x_{2}, x_{3},t ),t)\\
   R(x_{1}, x_{2}, x_{3},t)= \rho( x_{1}, x_{2}, \psi(x_{1}, x_{2}, x_{3},t),t)
   ,
  \end{split}
   \llabel{EQ11}
   \end{align}
for $x\in\Omega$ and~$t\geq0$.
The Euler equation in the new variables then become
  \begin{align}
   \partial_{t}v_{i}  
      +  v_{j} a_{kj} \partial_{k} v_{i} 
      - \psi_{t} a_{33} \partial_{3} v_{i}
      + \frac{1}{R} a_{ki}\partial_{k}  q(R) &= 0 
   \inin{\Omega}
   \comma i=1,2,3
   \label{EQ12}
   \\
   R_{t} + R a_{ji}\partial_{j}v_{i}+   v_{i} a_{ji}\partial_{j} R - \psi_{t} a_{33}\partial_{3} R& =0
   \inin{\Omega}
   ,
   \label{EQ13}
  \end{align}
where $a$ denotes the inverse of the Jacobian matrix~$\nabla \eta$.
We apply the summation convention on repeated indices throughout.
The function $q$ is given and satisfies
  \begin{equation}
    c_1 \leq   q'(R) \leq c_2
    \onon{[m_0/2,2M_0]}    
   .
   \label{EQ102}
  \end{equation}
We assume that $m_0>0$, $M_0$, $c_1>0$, and $c_2$ are constants, allowing all quantities to depend on them and on the smoothness properties of $q$
without mention.
In terms of $\psi$, the matrix $a$ reads
  \begin{align}
   a=
   \begin{pmatrix}
   1 & 0 & 0 \\
   0 & 1 & 0 \\
   \frac{-\partial_{1} \psi}{ \partial_{3} \psi}&  \frac{-\partial_{2} \psi}{ \partial_{3} \psi}& \frac{1}{ \partial_{3} \psi}
   \end{pmatrix} 
   .
  \llabel{EQ14}
  \end{align}
We also introduce the notation
  \begin{align}
   b=  \begin{pmatrix}
   \partial_{3} \psi & 0 & 0 \\
   0 &  \partial_{3} \psi& 0 \\
   -\partial_{1} \psi&  -\partial_{2} \psi& 1
   \end{pmatrix} 
   \llabel{EQ15}
  \end{align}
for the cofactor matrix corresponding to $a$,
while 
  \begin{align}
   J = \partial_{3} \psi
   \llabel{EQ16}
  \end{align}
stands for the determinant of the matrix~$\nabla \eta$.

Note that the plate equation is already defined in the reference
domain~$\Gamma_1$, except that now we use the new variable $q$ to
denote the pressure. Thus the equation reads
  \begin{align}
   w_{tt}  + \Deltah^{2} w - \Deltah w_{t}
    &= q(R)  \onon{\Gamma_1}
    ,
   \label{EQ17}
  \end{align}
assuming the periodic boundary conditions in the $x_{1}$ and $x_{2}$ directions.
The boundary conditions over the reference interface then become
 \begin{align}
   \label{EQ18} w_{t}    &= b_{3i} v_{i}   \onon{\Gamma_{1}}, \\
   \label{EQ19} v_{3}   &=  0              \onon{\Gamma_{0}}. 
\end{align}
The height $\psi $ satisfies the transport equation
  \begin{align}
   \begin{split}
   &
   \psi_{t}  = v_{3} - v_{2}\partial_{2} \psi-  v_{1}\partial_{1} \psi \onon{\Gamma_{1}}
   \\&
   \psi(\cdot,0) = 1
   \onon{\Gamma_1}
   .
  \end{split}
  \llabel{EQ20}
  \end{align}

\cole
\begin{Theorem}
\label{T01}
(A~priori estimates for existence)
Assume that
$(v,R,w)$ is a $C^{\infty}$ solution on an interval $[0,T]$ with
  \begin{align}
   \begin{split}
    \Vert v_0\Vert_{H^{3}},
    \Vert R_0\Vert_{H^{3}},
    \Vert R_0|_{\Gamma_1}\Vert_{H^{4}(\Gamma_1)},
    \Vert w_1\Vert_{H^{6}(\Gamma_1)}
    \leq M   
    ,
   \end{split}
   \label{EQ21}
  \end{align}
where $M\geq~$
and
  \begin{equation}
   m_0\leq R_0 \leq M_0
   ,
   \label{EQ162}
  \end{equation}
for some $m>0$.
Then there exists $T_0$ 
depending on $M$
such that 
$v$ and $w$ satisfy
  \begin{align} 
    \begin{split} 
     &\Vert \partial_{t}^{j}v\Vert_{H^{3-j}}, 
      \Vert \partial_{t}^{j}R\Vert_{H^{3-j}}, 
      \Vert w\Vert_{H^{5}(\Gamma_1)}, 
      \Vert w_t\Vert_{H^{4}(\Gamma_1)}, 
      \Vert w_{tt}\Vert_{H^{3}(\Gamma_1)}, 
      \Vert w_{ttt}\Vert_{H^{2}(\Gamma_1)}, 
      \Vert w_{tttt}\Vert_{L^{2}(\Gamma_1)}
  \leq K \comma t\in[0,T']
  ,
  \end{split}
   \label{EQ22}
  \end{align}
for $j=0,1,2,3$,
where $T'=\min\{T,T_0\}$,
with additionally
  \begin{equation}
   \int_{0}^{T_0}
    \Vert w_{tttt}\Vert_{H^{1}(\Gamma_1)}^2
    \,ds
    \leq
    K
    ,
   \llabel{EQ23}
  \end{equation}
while $a$, $b$, and $\psi$ satisfy
  \begin{equation}
    \Vert \partial_{t}^{j}\psi\Vert_{H^{5.5-j}}
    \Vert \partial_{t}^{j}a\Vert_{H^{4.5-j}},
    \Vert \partial_{t}^{j}b\Vert_{H^{4.5-j}}
    \leq
    K
    ,
   \llabel{EQ24}
  \end{equation}
for $j=0,1,2,3$,
where $K$ depends on~$M$.
\end{Theorem}
\colb

\startnewsection{Initial data}{sec3}
We consider initial data
as in \eqref{EQ21} such that
  \begin{equation}
   \partial_{t}R_{0}|_{\Gamma_{1}} \in H^{2}(\Gamma_{1})   
   ,
   \llabel{EQ27}
  \end{equation}
where $\partial_{t}R_{0}$ is defined below.
Consequently, 
we also have
  \begin{equation}
   (\psi(0),\psi_t(0)) \in H^{8.5}(\Omega) \times H^{6.5}(\Omega)
   .
   \llabel{EQ28}
  \end{equation}
We may then define iteratively the initial values of time derivatives as follows.
Starting with the continuity equation \eqref{EQ13} and evaluating at $t=0$, we define
  \begin{align}
   R_t(0)= - \div(R_{0}v_{0}) + \psi_{t}(0) \partial_{3}R_{0}
   ,
   \llabel{EQ29}
  \end{align}
from which we conclude $ R_t(0) \in H^2$.
From the
Euler equation \eqref{EQ12} at $t=0$, we set
  \begin{align}
   v_t(0)= -v_{0} \nabla v_{0}+ \psi_{t}(0) \partial_{3}v_{0} - \frac{1}{R_{0}}\nabla ( q(R_{0}))
   ,
   \llabel{EQ30}
  \end{align}
from which we obtain $v_t(0) \in H^2$. 
Taking the time derivatives repeatedly, we may similarly define the
initial time derivatives of $R$ and $v$ and obtain the regularities
  \begin{equation}
   \partial^{j}_{t} R(0) \in H^{3-j}
   \andand
   \partial^{j}_{t} v(0) \in H^{3-j}   
   ,
   \llabel{EQ31}
  \end{equation}
for $j= 0,1,2,3$.

Next, from the plate equation \eqref{EQ17}, we may define 
the initial time second derivative  as
  \begin{equation}
   w_{tt}(0)= -\Deltah^{2}w_{0}+ \Deltah w_{1} + q(R_{0})
   ,
   \llabel{EQ32}
  \end{equation}
leading to $w_{tt}(0) \in H^4$. 
Taking a time derivative of \eqref{EQ17}, we then set
  \begin{equation}
    w_{ttt}(0)= -\Deltah^{2} w_{1}+ \Deltah w_{tt} (0) + q'(R_{0})R_{t}(0)
    ,
   \llabel{EQ33}
  \end{equation}
from where $w_{ttt}(0) \in H^2$.
Similarly, we define
  \begin{equation}
   w_{tttt}(0)
   =
   -\Deltah^{2} w_{tt}(0)+ \Deltah w_{ttt}(0) + q''(R_{0})R^{2}_{t}(0)+ q'(R_{0})R_{tt}(0)
   ,
   \llabel{EQ34}
  \end{equation}
which gives $w_{tttt}(0) \in L^{2}$.

We may then introduce the total energy functional for the initial norms
  \begin{align}
   E(0)
   =
   \sum_{j=0}^{3}
      (
        \Vert  \partial^{j}_{t} v(0) \Vert^{2}_{H^{3-j}}
         +  \Vert  \partial^{j}_{t} R(0) \Vert^{2}_{H^{3-j}}
      )
   + \sum_{j=0}^{4} \Vert  \partial^{j}_{t} w(0) \Vert^{2}_{H^{8-2j}(\Gamma_{1})} 
   +\Vert R_{0} \Vert^{2}_{H^{4}(\Gamma_{1})} +\Vert  R_t(0) \Vert^{2}_{H^{2}(\Gamma_{1})}
   .
   \llabel{EQ35}
  \end{align}

\startnewsection{Preliminary estimates}{sec35}
The first statement addresses the continuity properties of $\psi$
and~$b$.

\cole
\begin{Lemma}
\label{L02}
(i) The harmonic extension $\psi$ defined in \eqref{EQ09} satisfies the estimates
  \begin{equation}
   \Vert \psi \Vert_{H^{s+1/2}} \lec \Vert w \Vert_{H^{s}(\Gamma_{1})}
       \comma s \in[0,5]
  \llabel{EQ145}
  \end{equation}
and
  \begin{equation}
   \Vert \partial^{j}_{t} \psi\Vert_{H^{s+1/2}}
   \lec
    \Vert \partial^{j}_{t} w \Vert_{H^{s}(\Gamma_{1})}
    \comma s\in[0,5-j]
    \commaone j=0,1,2,3
    .
   \llabel{EQ144}
   \end{equation}
(ii) The coefficient matrix $b$ satisfies 
  \begin{equation}
  \Vert \partial^{j}_{t} b \Vert_{H^{s-1/2}}
  \lec
  \Vert \partial^{j}_{t} w \Vert_{H^{s}(\Gamma_{1})}
    \comma s\in[1,5-j]
    \commaone j=0,1,2,3
  .
  \llabel{EQ146}
  \end{equation}
\end{Lemma}
\colb

The proof is obtained by a direct application of elliptic regularity estimates and it is thus omitted.

The next statement shows how to control $J$, $a$, and $R$ for small times
depending on the size of the solution $(v,R,w)$.

First, we introduce our polynomial notation.
Namely, in the rest of the paper, $\mathcal{P}$ denotes
a generic positive polynomial
of the variable
  \begin{align}
  \begin{split}
   &
    \sum_{j=0,1,2,3} \Bigl(    \Vert \partial_{t}^{j} v \Vert_{H^{3-j}}
    +
         \Vert \partial_{t}^{j} w \Vert_{H^{5-j}(\Gamma_1)}
    +
    \sup_{j=0,1,2,3}     \Vert \partial_{t}^{j} R \Vert_{H^{3-j}}
    \Bigr)
    + \Vert \partial_{t}^{4}w\Vert_{L^2(\Gamma_1)}
    ,
 \end{split}
   \llabel{EQ182}
  \end{align}
while $\mathcal{P}_0$ denotes a generic positive polynomial in
$M$ from~\eqref{EQ21}. Otherwise, we denote by $P$ any polynomial in the indicated arguments.

\cole
\begin{Lemma}
\label{L03}
Suppose that \eqref{EQ162} and \eqref{EQ22} hold on some interval
$[0,T]$.
Then
there exists $T'>0$ depending on $M$ and $T$ such that
the following
statements hold for $t\in [0,T']$:\\
(i) $1/2\leq J \leq 2$, 
\\
(ii)
$\Vert \partial^{j}_{t} a \Vert_{H^{s-1/2}} \lec \Vert \partial^{j}_{t} w \Vert_{H^{s}(\Gamma_{1})}P(\Vert \partial^{j-1}_{t} w \Vert_{H^{s}(\Gamma_{1})},\ldots, \Vert w \Vert_{H^{s}(\Gamma_{1})})$, for $j=0,1,2,3$ and $s\in[0,5-j]$,
\\
(iii) $m_0/2\leq R\leq 2M_0$,
\\
(iv) $    c_1 \leq   q'(R) \leq c_2$,
\\
(v)
$\Vert\partial_{t}^{j}(q(R))\Vert_{H^{3-j}}
\lec
\Vert\partial_{t}^j R\Vert_{H^{3-j}}
+ \mathcal{P}_0 + \int_{0}^{t} \mathcal{P}\,ds
$,
for $j=0,1,2,3$,  and
\\
(vi) for any $\epsilon\in(0,1]$ we have $\Vert a-I\Vert_{H^{2}} \leq \epsilon$
   provided $T'$ also depends on~$\epsilon$.
\end{Lemma}
\colb

\begin{proof}[Proof of Lemma~\ref{L03}]
(i) Note that $J = \partial_{3} \psi = 1+ \int_{0}^{t} \partial_{3} \psi_{t} $, so we have
  \begin{align}
  \begin{split}
  J  & \geq 1 -  \int_{0}^{t}  \Vert \nabla  \psi_{t} \Vert_{L^{\infty}}
	 \geq 1 -  C\int_{0}^{t}  \Vert \nabla  \psi_{t} \Vert_{H^{2}}
					\geq 1 -  C\int_{0}^{t}  \Vert w_{t} \Vert_{H^{5/2}(\Gamma_{1})}
				 \geq 1 - C t M
   ,
  \end{split}
  \llabel{EQ149}
  \end{align} 
and $(i)$ follows.

(ii) This follows from Lemma~\ref{L02}~(ii)
and the estimate on $J$ from the part (i).

(iii) We rewrite $R$ using the Fundamental Theorem of Calculus (FTC) as
  \begin{align} 
  \begin{split}
    R(t) & = R_0+\int_{0}^{t} R_{t}\,ds
    \geq m_0 - \int_{0}^{t}  \Vert   R_{t} \Vert_{L^{\infty}}
     \geq m_0 -  C \int_{0}^{t}  \Vert R_{t} \Vert_{H^{2}}
     \geq m_0 -  C t M
     ,
  \end{split}
  \llabel{EQ150}
  \end{align}
proving~(iii).

(iv) This follows from (iii) and~\eqref{EQ102}.

(v) This is a simple application of the chain rule; when differentiating
$\partial_t^{j} q(R)$, the leading term is
$q'(R)\partial_t^{j}R$, and all the lower order terms are estimated using the~FTC.

(vi) This follows by using the FTC in~$t$,
similarly to the proofs in (i) and~(ii).
\end{proof}

\startnewsection{a~priori estimates}{sec4}

The goal of this section to prove the following
time-tangential a~priori estimate.
It is convenient to introduce the notation
  \begin{equation}
   \bar\partial = (\partial_1,\partial_{2})
   \llabel{EQ118}
  \end{equation}
for the tangential gradient.

\cole
\begin{Lemma}
\label{L06}
Under the conditions of Theorem~\ref{T01}, we have
  \begin{align}
  \begin{split}
  &
    \sum_{j=0}^{3}
    \Bigl(
     \Vert \partial_{t}^{j}\bar\partial^{3-j} v(t) \Vert_{L^{2}}^{2}
     +   \Vert \partial_{t}^{j}\bar\partial^{3-j} R(t) \Vert_{L^{2}}^{2}
     + \Vert \partial^{j+1}_{t} \bar\partial^{3-j} w(t) \Vert_{L^{2}(\Gamma_1)}^{2} 
     +  \Vert \partial^{j}_{t} \bar\partial^{3-j} w\Vert_{H^{2}(\Gamma_1)}^{2} 
     \Bigr)
    \\&\indeq\indeq
   +
    \sum_{j=0}^{3}
    \int_{0}^{t} \Vert \partial^{j+1}_{t} \bar\partial^{4-j} w \Vert_{L^{2}(\Gamma_1)}^{2} 
    \\&\indeq
     \lec
    \epsilon \Vert q_{tt} \Vert_{H^{1}}^{9/5}
    +
     \mathcal{P}_0
     +\int_{0}^{t}\mathcal{P}\,ds
   .
  \end{split}
   \label{EQ56}
  \end{align}
\end{Lemma}
\colb

The proof is provided in Sections~\ref{sec51}--\ref{sec53}.

\subsection{Time derivative estimate}
\label{sec51}
Multiplying the momentum equation
\eqref{EQ12} with $JR$, we may rewrite it as
 \begin{equation}
   JR\partial_{t}v_{i}  +  Rv_{j} b_{kj} \partial_{k} v_{i} - R\psi_{t} b_{33} \partial_{3} v_{i}+  b_{ki}\partial_{k}  q = 0
   \inin{\Omega  \times (0,T)}
   \comma i=1,2,3
   .
   \llabel{EQ36}
  \end{equation}
We then apply the operator $\partial_{t}^{3}$ to the equation and take the $L^{2}$ inner product with $\partial_{t}^{3} v$ to obtain
  \begin{equation}
  \int\partial_{t}^{3}( JR  \partial_{t}v_i) \partial_{t}^{3} v_i
  = -\int \partial_{t}^{3}(  Rv_{j} b_{kj} \partial_{k} v_{i} )\partial_{t}^{3} v_i
  + \int\partial_{t}^{3}(R\psi_{t} b_{33} \partial_{3} v_{i} )\partial_{t}^{3} v_i -  \int \partial_{t}^{3}(b_{ki}\partial_{k}  q)\partial_{t}^{3} v_i
   ;
  \label{EQ37}
   \end{equation}
here and in the sequel, the integrals and norms without an indication of the
domain are understood to be over~$\Omega$.
We use the identity
  \begin{align}
  \begin{split}
    &    \frac{1}{2} \frac{d}{dt} \int JR  |\partial_{t}^{3} v|^{2} 
    = \frac12 \int \partial_{t}(JR)  |\partial_{t}^{3} v|^{2} 
      - \int
        \bigl(
             \partial_{t}^{3}(JR \partial_{t}v_{i}) 
             - JR \partial_{t}^{4} v_{i}
  \bigr)\partial_{t}^{3} v_{i} 
       + \int \partial_{t}^{3}( JR  \partial_{t}v_i) \partial_{t}^{3} v_i
  \end{split}
   \label{EQ38}
  \end{align}
and then replace the last term in \eqref{EQ38} with the expression in \eqref{EQ37} and rewrite the terms using commutators to get
  \begin{align}
  \begin{split}
    \frac{1}{2} \frac{d}{dt} \int JR  |\partial_{t}^{3} v|^{2} 
    & =   
      \underbrace{ \frac{1}{2}  \int 
         \partial_{t}(JR)  |\partial_{t}^{3} v|^{2} }_{K_{1}}
    - \underbrace{\int ( \partial_{t}^{3}(JR \partial_{t}v_{i}) 
        - JR \partial_{t}^{4} v_{i} )\partial_{t}^{3} v_{i} }_{K_{2}}
    \\ 
  & \indeq     
     -\underbrace{ \frac12\int R v_{j} b_{kj} \partial_{k} |\partial_{t}^{3} v|^{2}   + \frac{1}{2} \int R \psi_{t}  \partial_{3} |\partial_{t}^{3} v|^{2} }_{I_{1}}  
  \\
  & \indeq
   -    \underbrace{ \int (\partial_{t}^{3}
   (R v_{j} b_{kj} \partial_{k}v_{i}) - R v_{j} b_{kj} \partial_{k}\partial_{t}^{3}v_{i}  ) \partial_{t}^{3} v_{i} }_{K_{3}}
   \\
  & \indeq
   +   \underbrace{ \int (\partial_{t}^{3}
   (R \psi_{t}  \partial_{3}v_{i}) - R \psi_{t}  \partial_{3}\partial_{t}^{3}v_{i}  ) \partial_{t}^{3} v_{i}}_{K_{4}}
  \\
  & \indeq
  + \underbrace{\int b_{ki} \partial_{t}^{3} q \partial_{k}\partial_{t}^{3} v_{i}}_{I_{2}}
  - \underbrace{  \int (\partial_{t}^{3}(b_{ki}\partial_{k} q) - b_{ki} \partial_{t}^{3} \partial_{k}q) \partial_{t}^{3} v_{i}}_{K_{5}}
  -  \underbrace{\int_{\Gamma_{1}} b_{3i} \partial_{t}^{3} q \partial_{t}^{3} v_{i}}_{I_{\text{B}1}}
   ;
   \end{split}
   \label{EQ39} 
  \end{align}
we have also integrated by parts in $x_k$
to rewrite
$\int \partial_{t}^{3}(b_{ki}\partial_{k} q) \partial_{t}^{3} v_{i}$,
noting that one of the terms vanishes by the Piola identity $\partial_{k}b_{ki}=0$ for $i=1,2,3$.
The terms labeled as $K$ are of lower order in the sense that they satisfy the estimate
  \begin{align}
   \sum_{i=1}^{5} K_{i}
   \lec
   \mathcal{P}
   .
   \llabel{EQ40}
  \end{align}
The terms in $I_{1}$ are dealt with using integration by parts in space and observing that the boundary terms on $\Gamma_1$ vanish
due to \eqref{EQ18} and $\psi_t=w_t$ on~$\Gamma_1$; note that they also vanish on $\Gamma_0$ by~\eqref{EQ19} and $\psi_t=0$ on~$\Gamma_0$.
Using also the Sobolev and H\"older's inequalities, we may then estimate $I_{1}$ as
  \begin{align}
    I_{1}
    \lec
    \Vert \partial^{3}_{t} v \Vert^{2}_{L^{2}} \Vert R \Vert_{H^{3}} (  \Vert v \Vert_{H^{3}} \Vert b \Vert_{H^{2}}+\Vert \psi_{t} \Vert_{H^{3}}) 
    .
   \llabel{EQ41}
  \end{align} 
The high order terms $I_{2}$ and $I_{\text{B}1}$ 
are eventually canceled when combining with the density and plate estimates.
Summarizing the bounds on the terms in \eqref{EQ39}, we get
after integration in time
  \begin{align}
  \begin{split}
    &
     \frac12
     \int JR  |\partial_{t}^{3} v|^{2}
    \leq
    - \int_{0}^{t}I_{\text{B}1} \,ds
    + \int_{0}^{t}I_2 \,ds
    + \mathcal{P}_0
    + \int_{0}^{t}\mathcal{P} \,ds
    .
  \end{split}
   \label{EQ155}
  \end{align}

Performing the energy estimates on the plate equation \eqref{EQ17} after three time differentiations yields the identity
  \begin{align}
   \frac{1}{2} \frac{d}{dt} \Vert \partial^{4}_{t} w\Vert_{L^2(\Gamma_1)}^{2} +  \frac{1}{2} \frac{d}{dt} \Vert \Deltah \partial^{3}_{t} w\Vert_{L^2(\Gamma_1)}^{2} + \Vert \bar\partial \partial_{t}^{4} w \Vert_{L^2(\Gamma_1)}^{2} 
   = \underbrace{\int_{\Gamma_{1}} \partial_{t}^{3} q \partial_{t}^{4} w}_{I_{\text{B}}}
  .
  \label{EQ42}
  \end{align}
The boundary integral term $I_B$ on the right-hand side of \eqref{EQ42} may be expressed using the kinematic boundary condition \eqref{EQ18} as 
  \begin{align}
  I_{\text{B}}
  =
  \int_{\Gamma_{1}} \partial_{t}^{3} q \partial_{t}^{3} (b_{3i} v_{i}) 
  &=  \underbrace{\int_{\Gamma_{1}} \partial_{t}^{3} q b_{3i} \partial_{t}^{3}  v_{i}}_{I_{\text{B}1}}  
    + \underbrace{ \int_{\Gamma_{1}} \partial_{t}^{3} q    \partial_{t}^{3} b_{3i} v_{i}  }_{I_{\text{B}2}}
    +\underbrace{3  \int_{\Gamma_{1}} \partial_{t}^{3} q    \partial_{t} b_{3i}  \partial_{t}^{2} v_{i}  }_{I_{\text{B}3}}    
    +\underbrace{3  \int_{\Gamma_{1}} \partial_{t}^{3} q    \partial_{t}^{2} b_{3i}   \partial_{t} v_{i}  }_{I_{\text{B}4}}      
  .
  \label{EQ43}
  \end{align}
Note that these four terms cannot be estimated directly since they
  contain $\partial_{t}^{3}q$, which cannot be directly estimated on~$\Gamma_1$.
The first term $I_{\text{B}1}$ cancels with the first term in \eqref{EQ155} upon adding the equations \eqref{EQ155} and~\eqref{EQ42}.
We re-express $I_{\text{B}2}$ as
  \begin{equation}I_{\text{B}2}= \frac{d}{dt} \left(\int_{\Gamma_{1}} \partial_{t}^{2} q    \partial_{t}^{3} b_{3i} v_{i}\right)
        - \int_{\Gamma_{1}} \partial_{t}^{2} q    \partial_{t} (\partial_{t}^{3} b_{3i} v_{i})
   .
   \llabel{EQ44}
  \end{equation}
Integrating in time and employing H\"older's and Sobolev inequalities, we get
  \begin{align} 
  \begin{split}
   \int_{0}^{t} I_{\text{B}2} 
     &=  
       \int_{\Gamma_{1}} \partial_{t}^{2} q(t)  \partial_{t}^{3} b_{3i}(t) v_{i}(t) 
       - \int_{\Gamma_{1}} \partial_{t}^{2} q(0)  \partial_{t}^{3} b_{3i}(0) v_{i}(0)
       - \int_{0}^{t} \int_{\Gamma_{1}}  \partial_{t}^{2} q   \partial_{t} (\partial_{t}^{3}b_{3i} v_{i})
    .
  \end{split}
   \label{EQ45}
  \end{align}
The first term on the right may be estimated as
  \begin{align}
  \begin{split}
   \int_{\Gamma_{1}} \partial_{t}^{2} q(t)  \partial_{t}^{3} b_{3i}(t) v_{i}(t) 
   &
   \lec
   \Vert \partial_{t}^2 q\Vert_{H^{1/2}(\Gamma_1)}
   \Vert \partial_{t}^3 b_{3i}\Vert_{L^2(\Gamma_1)}
   \Vert v_i\Vert_{H^{1/2}(\Gamma_1)}
   \\&
   \lec
   \Vert \partial_{t}^2 q\Vert_{H^{1}}
   \Vert \partial_{t}^3 \psi\Vert_{H^2}^{2/3}
   \Vert \partial_{t}^3 \psi\Vert_{H^{1/2}}^{1/3}
   \Vert v\Vert_{H^{1}}
   \\&
   \lec
   \epsilon    \Vert \partial_{t}^2 q\Vert_{H^{1}}^{9/5}
   +
   \epsilon    \Vert \partial_{t}^3 \psi\Vert_{H^2}^{2}
   + C_{\epsilon}\mathcal{P}_0
   + C_{\epsilon}\int_{0}^{t}\mathcal{P}\,ds
   ,
  \end{split}
   \llabel{EQ113}
  \end{align}
for any $\epsilon\in(0,1]$;
we used $b_{3i} = [ -\partial_{1} \psi, -\partial_{2} \psi, 1]$, which contains
only tangential derivatives,
and thus we  estimated
$
\Vert \partial_{t}^{3}b_{3i}\Vert_{L^{2}(\Gamma_1)}
\lec
\Vert \partial_{t}^{3}\bar\partial \psi\Vert_{L^{2}(\Gamma_1)}
 \lec
 \Vert \partial_{t}^{3}\psi\Vert_{H^{3/2}} 
$.
We used the power $9/5$ since a power less than 2 is needed in~\eqref{EQ114}.
In summary, we get
  \begin{align} 
  \begin{split}
   \int_{0}^{t} I_{\text{B}2} 
    &
    \lec
       \epsilon \Vert  \partial_{t}^{2} q(t) \Vert_{H^{1}}^{9/5}
       + \epsilon \Vert \partial_{t}^{3}  \psi \Vert^{2}_{H^{2}}
       + \epsilon  \int_{0}^{t} \Vert \partial_{t}^{4}  \psi \Vert^{2}_{H^{1}(\Gamma_{1})}
       + C_{\epsilon}\mathcal{P}_0
       + C_{\epsilon}\int_{0}^{t} \mathcal{P}
    ,
  \end{split}
   \label{EQ114}
  \end{align}
for any $\epsilon\in(0,1]$.
Note that the third term on the right-hand side results from the last term in~\eqref{EQ45}. Namely,
  \begin{align}
  \begin{split}
   &
      - \int_{0}^{t} \int_{\Gamma_{1}}  \partial_{t}^{2} q\partial_{t}^{4}b_{3i} v_{i}
    \lec
    \int_{0}^{t}
    \Vert \partial_{t}^2 q\Vert_{L^2(\Gamma_1)}
    \Vert \partial_{t}^{4}b_{3i}\Vert_{L^2(\Gamma_1)}
    \Vert v\Vert_{L^{\infty}}
    \,ds
    \\&\indeq
    \lec
    \int_{0}^{t}
    \Vert \partial_{t}^2 q\Vert_{H^{1}}
    \Vert \bar\partial\partial_{t}^{4} \psi\Vert_{L^2(\Gamma_1)}
    \Vert v\Vert_{H^{2}}
    \,ds
    \lec
       \epsilon  \int_{0}^{t} \Vert \bar\partial\partial_{t}^{4}  \psi \Vert^{2}_{L^2(\Gamma_{1})}
       + C_{\epsilon}\int_{0}^{t} \mathcal{P}
    \\&\indeq
     \lec
       \epsilon  \int_{0}^{t} \Vert \partial_{t}^{4}  \bar \partial w \Vert^{2}_{L^2(\Gamma_{1})}
       + C_{\epsilon}\int_{0}^{t} \mathcal{P}
    ,
    \end{split}
   \label{EQ165}
  \end{align}
while when the time derivative in the last term of \eqref{EQ45} is applied to $v_i$, we may bound the term simply by
$\int_{0}^{t} \mathcal{P}$.

To address the term $I_{\text{B}3}$, we need to apply the divergence theorem to obtain
the interior integral and then integrate by parts in time in high order terms.
Namely, let $\zeta=\zeta(x_3)$ be a smooth cut-off function, which equals $1$ in a neighborhood of $x_3=1$ and
vanishes for $x_3\leq 1/2$.
Then we have
  \begin{align}
  \begin{split}
    \frac13
    \int_{0}^{t}
    I_{\text{B}3}
    \,ds
    &=
    \int_{0}^{t}     
      \int \partial_{t}^{3} \partial_{3}q    \partial_{t} b_{3i}  \partial_{t}^{2} v_{i} \zeta
    +     \int_{0}^{t}      \int \partial_{t}^{3}q    \partial_{t}  \partial_{3}b_{3i}  \partial_{t}^{2} v_{i} \zeta
    \\&\indeq
    +     \int_{0}^{t}      \int \partial_{t}^{3}q    \partial_{t} b_{3i}  \partial_{t}^{2}  \partial_{3}v_{i} \zeta
    +      \int_{0}^{t}     \int \partial_{t}^{3}q    \partial_{t} b_{3i}  \partial_{t}^{2} v_{i}    \zeta'(x_3)
   \\&
   \leq
   \int \partial_{t}^{2} \partial_{3}q    \partial_{t} b_{3i}  \partial_{t}^{2} v_{i} \zeta \Big|_{0}^{t}
   + \int_{0}^{t}\mathcal{P}\,ds
   ,
   \end{split}
   \llabel{EQ147}
  \end{align}
where we integrated in time in the first term.
Therefore,
  \begin{align}
  \begin{split}
   \int_{0}^{t}
   I_{\text{B}3}
   \,ds
   &\lec
   \mathcal{P}_0
   + \Vert \partial_{t}^{2}q\Vert_{H^{1}}
     \Vert \partial_{t}b\Vert_{H^{2}}
     \Vert \partial_{t}^2 v\Vert_{L^2}
   +  \int_{0}^{t}\mathcal{P}\,ds
   \\&
   \lec
   \mathcal{P}_0
   + \epsilon \Vert \partial_{t}^{2}q\Vert_{H^{1}}^{9/5}
   + C_\epsilon \mathcal{P}_0
   + C_\epsilon \int_{0}^{t}\mathcal{P}\,ds
   .
  \end{split}
   \label{EQ187}
  \end{align}

The term $I_{\text{B}4}$ is handled analogously.
Here, the pointwise in time term reads
$   \int \partial_{t}^{2} \partial_{3}q    \partial_{t}^2 b_{3i}  \partial_{t} v_{i} \zeta $, which we bound as
  \begin{align}
  \begin{split}
   \int
      \partial_{t}^{2}\partial_{3}q
      \partial_{t}^2 b_{3i}
      \partial_{t} v_{i} \zeta
   \lec
   \Vert \partial_{t}^{2}q\Vert_{H^{1}}
   \Vert \partial_{t}^2 b\Vert_{H^{1}}
   \Vert \partial_{t} v\Vert_{H^{1}}
   ,
  \end{split}
   \llabel{EQ115}
  \end{align}
leading to
the same estimate, i.e.,
  \begin{align} 
  \begin{split}
   \int_{0}^{t} I_{\text{B}4}  \,ds
      &\lec 
   \mathcal{P}_0
   + \epsilon \Vert \partial_{t}^{2}q\Vert_{H^{1}}^{9/5}
   + C_\epsilon \mathcal{P}_0
   + C_\epsilon \int_{0}^{t}\mathcal{P}\,ds
   . 
  \end{split}
   \label{EQ46}
  \end{align}
Integrating \eqref{EQ42} and summarizing the bounds on the right-hand side, 
we get
  \begin{align}
    \begin{split}
    &
    \frac12
    \Vert \partial^{4}_{t} w\Vert_{L^2(\Gamma_1)}^{2}
     +
    \frac12
     \Vert \Deltah \partial^{3}_{t} w\Vert_{L^2(\Gamma_1)}^{2}
     + \int_{0}^{t}\Vert \bar\partial \partial_{t}^{4} w \Vert_{L^2(\Gamma_1)}^{2} 
    \\&\indeq
    \leq
     I_{\text{B}1}
     +
       \epsilon    \Vert \partial_{t}^2 q\Vert_{H^{1}}^{9/5}
   +
   \epsilon    \Vert \partial_{t}^3 \psi\Vert_{H^2}^{2}
   + C_{\epsilon}\mathcal{P}_0
   + C_{\epsilon}\int_{0}^{t}\mathcal{P}\,ds
  ,
  \end{split}
   \label{EQ148}
  \end{align}
where we absorbed the first term on the far right side of \eqref{EQ165}
using the third term on the left-hand side of~\eqref{EQ42}.

Next, we apply $\partial_{t}^{3}$ to the continuity equation \eqref{EQ13} 
\begin{equation}
   J R_{t} + R b_{ji}\partial_{j}v_{i}+   v_{i} b_{ji}\partial_{j} R - \psi_{t} b_{33}\partial_{3} R =0
   \inin{\Omega}   
   \llabel{EQ49}
  \end{equation}
and integrate against $\frac{q'(R)}{R}\partial_{t}^{3}R$ to get an expression for
$ \int  \frac{q'(R)}{R}\partial_{t}^{3}\left( J  \partial_{t}R\right) \partial_{t}^{3} R$
which we use to replace the last term in the identity
  \begin{align}
  \begin{split}
    \frac{1}{2} \frac{d}{dt} \int J \frac{q'(R)}{R}  |\partial_{t}^{3} R|^{2}
    &= \frac12 \int \partial_{t}\left(J \frac{q'(R)}{R}\right)  |\partial_{t}^{3} R|^{2} 
      - \int \left( \partial_{t}^{3}\left(J \frac{q'(R)}{R} \partial_{t}R\right) 
                      - J \frac{q'(R)}{R} \partial_{t}^{4} R 
                      \right)\partial_{t}^{3} R
    \\&\indeq
       + \int \left(\partial_{t}^{3}
                     \left( J \frac{q'(R)}{R}  \partial_{t}R\right)
                     -
                     \frac{q'(R)}{R}
                  \partial_{t}^{3}  ( J \partial_{t}R)
                       \right)
             \partial_{t}^{3} R
       + \int  \frac{q'(R)}{R}\partial_{t}^{3}\left( J  \partial_{t}R\right) \partial_{t}^{3} R
   .
  \end{split}
   \llabel{EQ50}
  \end{align}
Rewriting the identity using commutators, we have
  \begin{align}
  \begin{split}
     &
    \frac{1}{2}\frac{d}{dt}\int     J \frac{q'(R)}{R}| \partial_{t}^{3}R |^{2}
    \\&
    =
      \underbrace{\frac12 \int \partial_{t}\left(J \frac{q'(R)}{R}\right)  |\partial_{t}^{3} R|^{2} }_{I_0}
      - 
       \underbrace{\int \left( \partial_{t}^{3}\left(J \frac{q'(R)}{R} \partial_{t}R\right) 
                      - J \frac{q'(R)}{R} \partial_{t}^{4} R 
                      \right)\partial_{t}^{3} R}_{K_6}
    \\&\indeq
       + \underbrace{\int \left(\partial_{t}^{3}
                     \left( J \frac{q'(R)}{R}  \partial_{t}R\right)
                     -
                     \frac{q'(R)}{R}
                  \partial_{t}^{3}  ( J \partial_{t}R)
                       \right)
             \partial_{t}^{3} R}_{K_7}
   - \underbrace{\int  b_{ji} \partial_{t}^{3} \partial_{j}v_{i}   \partial_{t}^{3}q(R)}_{I_{2}}
   \\& \indeq
   - \underbrace{ \frac{1}{2} \int  v_{i} b_{ji} \partial_{j}(\partial_{t}^{3}R)^{2}   \frac{q'(R)}{R}
   + \frac{1}{2}\int  \psi_{t}    \partial_{3}(\partial_{t}^{3}R)^{2}   \frac{q'(R)}{R} }_{I_{3}}
    \\&\indeq
    -  \underbrace{\int  \frac{q'(R)}{R}
       \Bigl(
        \partial_{t}^{3}(R b_{ji}\partial_{j}v_i)
        - R b_{ji}\partial_{t}^{3}\partial_{j}v_i
       \Bigr)
       \partial_{t}^{3} R}_{K_8}
    -  \underbrace{\int  
       b_{ji}\partial_{t}^{3}\partial_{j} v_i
       \Bigl(
         q'(R) \partial_{t}^{3} R - \partial_{t}^{3}(q(R))
       \Bigr)}_{K_9}
    \\&\indeq
   - \underbrace{\int
      \frac{q'(R)}{R} \Bigl(
        \partial_{t}^{3}  (v_i b_{ji}\partial_{j} R)
        - v_i b_{ji}\partial_{t}^{3}\partial_{j}R
       \Bigr)
       \partial_{t}^{3} R }_{K{10}}
   +
   \underbrace{
   \int
    \frac{q'(R)}{R}
    \Bigl(
     \partial_{t}^{3}(\psi_t b_{33}\partial_{3}R)
     - \psi_t b_{33} \partial_{t}^{3}\partial_{3}R
    \Bigr)
    \partial_{t}^{3} R}_{K_{11}}
   .
  \end{split}
   \label{EQ51}
  \end{align}
The terms $K_6$--$K_{11}$ are lower order commutator terms satisfying
  \begin{align}
   \sum_{m=6}^{11}\int_{0}^{t}  K_m
   \lec 
      \int_{0}^{t}\mathcal{P}
   .
   \label{EQ52}
  \end{align}
The term $I_{0}$ satisfies a similar estimate, while the terms in $I_{3}$ are estimated after integrating by parts in the space variable and noting that the boundary terms vanish due to the kinematic boundary condition. We thus obtain
  \begin{align}
   \int_{0}^{t}
   I_{3}
   \lec
   \int_{0}^{t}
   \mathcal{P}
   .
   \label{EQ53}
  \end{align}
Integrating \eqref{EQ51} in time and 
using \eqref{EQ52}--\eqref{EQ53}, 
we get
  \begin{align}
  \begin{split}
   &
   \frac{1}{2}\int     J \frac{q'(R)}{R}| \partial_{t}^{3}R |^{2}
    \leq
     - \int_{0}^{t} I_{2}
     + \mathcal{P}_0
     + \int_{0}^{t}\mathcal{P}
   .
  \end{split}
   \label{EQ54}   
  \end{align}
Adding the equations \eqref{EQ155}, \eqref{EQ148}, and~\eqref{EQ54}, 
while noting that the terms involving $I_{2}$ and $I_{\text{B}}$ cancel, we have
  \begin{align}
  \begin{split}
      &\Vert \partial_{t}^{3} v(t) \Vert_{L^{2}}^{2}
    +   \Vert \partial_{t}^{3} R(t) \Vert_{L^{2}}^{2}
    + \Vert \partial^{4}_{t} w(t) \Vert_{L^{2}(\Gamma_1)}^{2} +  \Vert \Deltah \partial^{3}_{t} w\Vert_{L^{2}(\Gamma_1)}^{2} 
    + \int_{0}^{t} \Vert \bar\partial \partial_{t}^{4} w \Vert_{L^{2}(\Gamma_1)}^{2} 
   \\&\indeq
   \lec
     \epsilon \Vert  \partial_{t}^{2} q(t) \Vert_{H^{1}}^{9/5}
    + \epsilon \Vert \partial_{t}^{3}  \psi \Vert^{2}_{H^{2}}
   +    C_\epsilon\mathcal{P}_0
    + C_\epsilon\int_{0}^{t}\mathcal{P}\,ds
   .
  \end{split}
   \label{EQ55}
  \end{align}
We also used Lemma~\ref{L03}~(iv) and~\eqref{EQ102}.
The last term on the right is bounded using
  \begin{align}
  \begin{split}
    \Vert \partial_{t}^{3}  \psi \Vert^{2}_{H^{3/2}}
    &\lec
    \Vert \partial_{t}^{3}  w \Vert^{2}_{H^{2}(\Gamma_1)}
    \lec
    \Vert \partial_{t}^{3}  \Deltah w \Vert^{2}_{L^{2}(\Gamma_1)}
    +
    \Vert \partial_{t}^{3}   w \Vert^{2}_{L^{2}(\Gamma_1)}    
    \lec
    \Vert \partial_{t}^{3}  \Deltah w \Vert^{2}_{L^{2}(\Gamma_1)}
    +
    \mathcal{P}_0
    + \int_{0}^{t}\mathcal{P}\,ds
    ,
    \end{split}
   \label{EQ77}
  \end{align}
and the first term on the far right may be absorbed into the left-hand side of \eqref{EQ55} once we multiply \eqref{EQ77} by~$\epsilon$.
Note that due to the second inequality in \ref{EQ77}, the strength of the damping in \eqref{EQ17} can be reduced; however, we do not pursue this direction further.
Using \eqref{EQ77} in \eqref{EQ55}, we finally obtain
  \begin{align}
  \begin{split}
     &\Vert \partial_{t}^{3} v(t) \Vert_{L^{2}}^{2}
    +   \Vert \partial_{t}^{3} R(t) \Vert_{L^{2}}^{2}
    + \Vert \partial^{4}_{t} w(t) \Vert_{L^{2}(\Gamma_1)}^{2} +  \Vert \Deltah \partial^{3}_{t} w\Vert_{L^{2}(\Gamma_1)}^{2} 
    + \int_{0}^{t} \Vert \bar\partial \partial_{t}^{4} w \Vert_{L^{2}(\Gamma_1)}^{2} 
   \\&\indeq
   \lec
    \epsilon \Vert  \partial_{t}^{2} q(t) \Vert_{H^{1}}^{9/5}
    +
      C_\epsilon\mathcal{P}_0
    +  C_\epsilon\int_{0}^{t}\mathcal{P}
   ,
  \end{split}
   \llabel{EQ61}
  \end{align}
completing the tangential energy estimate for the pure time derivatives.
\colb

\begin{proof}[Proof of Lemma~\ref{L06}]
The approach above can be repeated more generally replacing time derivatives with tangential derivatives in the $x_{1}, x_{2}$ directions.
In particular, we apply $\bar\partial^{l} \partial^{3-l}_{t}$ to each equation and use multipliers  $\bar\partial^{l} \partial^{3-l}_{t} ( v, w_{t}, R)$ for $s=1,2$ and $l=1,2,3$ respectively.
Note that we do not need to integrate by parts in time at any point and use
integration by parts in the tangential variables instead.
\end{proof}

\subsection{Alternative plate energy identity} 
\label{sec53}
The following statement is needed when combining with the hyperbolic
estimate on the density~$R$; see
the term
$I_{\BM4\text{b}}$ in \eqref{EQ108} below.

\cole
\begin{Lemma}
The function $w$ satisfies the inequality
  \begin{align} 
  \begin{split}
  & \frac1C \Vert \partial^{4}_{t} w \Vert^{2}
    + \frac1C \Vert \Deltah \partial^{3}_{t} w\Vert^{2}
    + \frac1C \int_{0}^{t}\Vert \bar\partial \partial_{t}^{4} w \Vert^{2} 
  \leq
  \underbrace{\int_{0}^{t}\int_{\Gamma_{1}} \frac{1}{JR} \partial_{t}^{3} q \partial_{t}^{4} w}_{I_{\text{B}}} 
  +  \mathcal{P}_0
  + \int_{0}^{t}\mathcal{P}\,ds
  ,
  \end{split}
  \label{EQ57}
  \end{align}
for some constant~$C\geq1$.
\end{Lemma}
\colb

\begin{proof}
Applying $\partial_{t}^{3}$ to the plate equation \eqref{EQ17} and 
integrating against $\frac{1}{JR}\partial_{t}^{4} w $ over $\Gamma_{1}$, we get
  \begin{align} 
  \begin{split}
  &\frac{1}{2} \frac{d}{dt} \left\Vert \frac{1}{\sqrt{JR}}\partial^{4}_{t} w \right\Vert_{L^2(\Gamma_1)}^{2} +  \frac{1}{2} \frac{d}{dt}  \left\Vert \frac{1}{\sqrt{JR}} \Deltah \partial^{3}_{t} w\right\Vert_{L^2(\Gamma_1)}^{2} + \left\Vert\frac{1}{\sqrt{JR}} \bar\partial \partial_{t}^{4} w \right\Vert_{L^2(\Gamma_1)}^{2} 
  \\&\indeq
  =
  - \int_{\Gamma_{1}} \frac{\partial_{t}(JR)}{J^{2}R^{2}}( |\partial^{4}_{t} w|^{2} +|\Deltah\partial^{3}_{t}   w|^{2})
  +\underbrace{\int_{\Gamma_{1}} \frac{1}{JR} \partial_{t}^{3} q \partial_{t}^{4} w}_{I_{\text{B}}} 
  \\&\indeq\indeq
  + \int_{\Gamma_{1}} 2 \frac{\nabla (JR)}{J^{2}R^{2}} \cdot
  ( \Deltah\partial^{3}_{t} w\,  \bar\partial \partial^{4}_{t} w
  + \partial^{4}_{t} w \,\bar\partial \partial^{4}_{t} w)
  + \int_{\Gamma_{1}} \div\left(\frac{\nabla (JR)}{J^{2}R^{2}}\right)  \partial^{4}_{t} w \Deltah\partial^{3}_{t} w
  ,
  \end{split}
   \llabel{EQ186}
  \end{align}
where the dot product is only in the horizontal directions.
Thus \refer{EQ57} follows by integrating in time and
using Lemma~\ref{L03}~(i) and~(iii).
\end{proof}
\colb

\startnewsection{A hyperbolic equation satisfied by the density}{sec5}

Denote by 
 \begin{equation}
   \nu=[b_{31},b_{32},b_{33}]^{T}
   ,
   \label{EQ58}
  \end{equation}
the moving normal. Also, we use the derivative notation
  \begin{equation}
   \div_{a} z=  a_{ki} \partial_{k} z_{i}
   \andand
   \nabla_{a} y = a_{km} \partial_{k} y
   ,
   \label{EQ59}
  \end{equation}
representing the variable divergence and gradient;
in \eqref{EQ59}, $z$ is a vector and $y$ a scalar function.
Introduce the differential operator 
  \begin{align}
   Q 
   =
   \partial_{t} + (v- \eta_{t}) \cdot \nabla_{a}
   = 
   \partial_{t} + (v_i- \partial_{t}\eta_{i}) a_{ki}\partial_{k}
   =
   \partial_{t} + v_i a_{ki}\partial_{k}
   - \psi_t a_{33}\partial_{3}
   ,
   \llabel{EQ60}
  \end{align}
where 
$\eta$ is as in \eqref{EQ10},
which is a first order differential operator of transport type.
Since the expression $v-\eta_t$ appears frequently, we introduce
the notation
  \begin{equation}
   \vv = v-\eta_t
   .
   \llabel{EQ135}
  \end{equation}
Moreover, the restriction of $Q$ to the boundary is tangential
due to the boundary conditions~\eqref{EQ18} and~\eqref{EQ19}. Namely,
  \begin{align}
   Q |_{\partial \Omega}
   =
   \partial_{t} + \sum_{k=1}^{2}v_k \partial_{k}
   .
   \llabel{EQ60b}
  \end{align}

In the next theorem, we derive a hyperbolic equation with Neumann-type boundary conditions 
for the logarithm of the density function~$R$.

\cole
\begin{Theorem}
\label{T02}
The function $g = \log R$ 
satisfies the hyperbolic equation
 \begin{align}
   Q^{2}g  
   - \div_{a} ( f \nabla_{a} g)    = F    \inin{\Omega  \times (0,T)}
  , \label{EQ62}
  \end{align}
with the boundary conditions
 \begin{align}
  \begin{split}
    f\nabla_{a}   g \cdot \nu = h  ~~\mbox{on}~~ {\Gamma_{1}  \times (0,T)} \\
    \partial_{3} g = 0  ~~\mbox{on}~~ {\Gamma_{0} \times (0,T)}
    ,
  \end{split}
   \label{EQ63}
  \end{align}
where $f= q'(R)$ and
  \begin{align}
   F=
     - \partial_{t}a_{ji} \partial_{j} v_{i}  
     - a_{mi} \partial_{m} (v_{j} a_{kj}) \partial_{k} v_{i}
     + \partial_{k} a_{mi} v_ja_{kj} \partial_{m}v_i
     +  a_{mi} \partial_{m}(\psi_t a_{33}) \partial_{3}v_i
     - \partial_{3}a_{mi} \psi_t a_{33}\partial_{m}v_i
   ,
   \label{EQ64}
  \end{align}
with
  \begin{align}
   \label{EQ65a}
   h =  \partial_{t}b_{3i} v_{i} -w_{tt}
   +  \sum_{j=1}^{2} \partial_{j} b_{3i} v_{l} a_{jl}  v_{i} 
   - \sum_{j=1}^{2} \partial_{j}w_{t} v_{l} a_{jl}  .
  \end{align}
\end{Theorem}
\colb

It is essential in the proof below that $h$ can be extended to $\Omega$,
and this is since we need to apply $Q$ and $Q^2$ to~$h$ in the interior.
We extend $h$ by introducing
  \begin{align}
   \label{EQ65}
   h = \partial_{t}b_{3i} v_{i}
   -  \psi_{tt}
   +  \sum_{j=1}^{2} \partial_{j} b_{3i} v_{l} a_{jl}  v_{i} 
   - \sum_{j=1}^{2} \partial_{j}\psi_{t} v_{l} a_{jl}
   ,
  \end{align}
noting that
$\psi_{t}|_{\Gamma_1}=w_{t}|_{\Gamma_1}$ and
$\psi_{tt}|_{\Gamma_1}=w_{tt}|_{\Gamma_1}$.

\begin{proof}
To obtain the boundary conditions for $R$, we take
the dot product of the momentum equation  \eqref{EQ12} with the moving normal vector 
\eqref{EQ58}
obtaining
  \begin{align}
  \begin{split}
   \frac{q'(R)}{R} b_{3i} a_{ki} \partial_{k}R 
    = -b_{3i} \partial_{t}v_{i} 
    -  b_{3i}v_{l} a_{jl} \partial_{j} v_{i} 
    + b_{3i} \psi_{t}a_{33} \partial_{3} v_{i}  
   .
   \end{split}
   \llabel{EQ66}
  \end{align}
Restricting to $\Gamma_{1}$, using the boundary conditions \eqref{EQ18}, and rewriting, we get
  \begin{align}
  \begin{split}
   \frac{q'(R)}{R} b_{3i} a_{ki} \partial_{k}R 
    &=
    \partial_{t}b_{3i} v_{i}
    - w_{tt}
    +  \partial_{j} b_{3i} v_{l} a_{jl}  v_{i} 
    - \partial_{3}b_{3i} \psi_{t}a_{33}  v_{i}  
    \\&\indeq
    -  \partial_{j}( b_{3i}v_{i} ) v_{l} a_{jl}   
    + \partial_{3}(b_{3i}v_{i}) \psi_{t}a_{33}    
   \onon{\Gamma_1}
   .
  \end{split}
  \label{EQ67}
  \end{align}
Note that when $j=3$, the
fifth term on the right cancels with the last term,
as
$\psi_{t}a_{33} =v_{l} a_{3l}$ by~\eqref{EQ18}.
Similarly, when $j=3$, the third term cancels with the fourth
due to
$v_l a_{3l}=\psi_t a_{33}$ on~$\Gamma_1$. Also, when
$j=1,2$ in the third term, we may use
$\partial_{j}( b_{3i}v_{i}) = \partial_{j}w_{t}$
on~$\Gamma_1$.
Using all these observations in \eqref{EQ67}, we obtain
  \begin{align}
  \frac{q'(R)}{R} b_{3i} a_{ki} \partial_{k}R 
  =  \partial_{t}b_{3i} v_{i} -w_{tt}
   +  \sum_{j=1}^{2} \partial_{j} b_{3i} v_{l} a_{jl}  v_{i} 
   - \sum_{j=1}^{2} \partial_{j}w_{t} v_{l} a_{jl}
   = h \onon{\Gamma_{1}}
  , \label{EQ68}
  \end{align}
and a similar computation on the rigid boundary   $\Gamma_{0}$ using \eqref{EQ19} where $v_{3}=0$ and $b_{3i} = \delta_{3i}$ and $\psi=0$, gives
  \begin{align}
   \partial_{3}R =  0    \onon{\Gamma_{0}}
   .
   \label{EQ69}
  \end{align}
  Using the notation \eqref{EQ58} and \eqref{EQ59}, we may express the 
Neumann boundary conditions \eqref{EQ68} and \eqref{EQ69} satisfied by $g$ as
  \begin{align}
  \begin{split}
    &
    q'(R) b_{3i} a_{ki} \partial_{k}g =
    f\nabla_{a}   g \cdot \nu = h
    \onon{\Gamma_{1}}
    \\&
    f\nabla_{a}   g \cdot \nu = f \partial_{3} g =0
    \onon{\Gamma_{0}}
    .
  \end{split}
   \llabel{EQ70}
  \end{align}

Next, we show that
  \begin{equation}
   g \equiv \log R   
   \llabel{EQ164}
  \end{equation}
satisfies a hyperbolic-type equation. First, denoting by
  \begin{align}
   H = a_{ki} \partial_{k} v_{i}
   ,
   \llabel{EQ71}
  \end{align}
the ALE divergence of $v$,
we may rewrite the continuity equation \eqref{EQ13}  in terms of $g$ as
  \begin{align}
   \llabel{EQ72}
   \partial_{t} g + v_{i} a_{ki} \partial_{k} g - \psi_{t}a_{33} \partial_{3} g + H =0
   .
  \end{align}
This may be expressed using a vector notation as 
  \begin{equation}
     \partial_{t} g + \vv\cdot \nabla_{a} g  + H =0   ,   
   \llabel{EQ73}
  \end{equation}
or simply
  \begin{align}
   \label{EQ74}
   Qg + H=0
   .
  \end{align}
On the other hand, we may rewrite the momentum equations \eqref{EQ12} as
  \begin{align}
   \partial_{t} v_{i} + v_{j} a_{kj} \partial_{k} v_{i}- \psi_{t}a_{33} \partial_{3} v_{i}   + a_{ki} q'(R) \partial_{k} g =0
   \comma i=1,2,3
   .
   \label{EQ75}
   \end{align}  
Applying $a_{ki}\partial_{k}$
to \eqref{EQ75}, we obtain 
the equation
  \begin{align}
  \begin{split}
   &
   \partial_{t} H 
     + v_{j} a_{kj} \partial_{k} H 
     -\psi_{t}a_{33} \partial_{3} H 
     + a_{ji} \partial_{j} (a_{ki} q'(R) \partial_{k} g) 
   = -F
   ,
  \end{split}
   \llabel{EQ76a}
  \end{align}  
where $F$ is given in \eqref{EQ64},
which may be, by using the definition of $Q$, expressed as 
  \begin{align}
   QH  
    + a_{ji} \partial_{j} (a_{ki} q'(R) \partial_{k} g) 
    =- F
   .
   \label{EQ78}
 \end{align}  
Applying $Q$ to the equation \eqref{EQ74} and substituting the expression for $QH$ from \eqref{EQ78} while denoting $f=q'(R)$, we get~\eqref{EQ62}.
\end{proof}

\cole
\begin{Lemma}
\label{hF}
The function $F$ defined in \eqref{EQ64} and the trace $h$ defined in \eqref{EQ65} satisfy the estimate
  \begin{align}
  \begin{split}
  \sum_{k=1}^{2}
  \Vert \partial_{t}^{2-k} h \Vert_{H^{0.5+k}(\Gamma_{1})}
  +
  \sum_{k=0}^{2}\Vert \partial_{t}^{2-k} F \Vert_{H^{k}}
  \lec
    \mathcal{P}
  ,
  \end{split}
   \label{EQ151}
  \end{align}
while
  \begin{align}
  \begin{split}
  \Vert \partial_{t}^{2} h \Vert_{H^{0.5}(\Gamma_{1})} 
  \lec
   \Vert \partial_{t}^{4} w \Vert_{H^{0.5}(\Gamma_{1})} 
   +
    \mathcal{P}
  ,
  \end{split}
   \label{EQ166}
  \end{align}
\end{Lemma}
\colb

In \eqref{EQ151} and \eqref{EQ166}, as well as in many inequalities
below,
the quantities are understood to be evaluated at a time $t\in[0,T]$.

\begin{proof}
Both inequalities \eqref{EQ151} and \eqref{EQ166} are obtained
directly 
using H\"older and Sobolev inequalities. We omit further details.
\end{proof}

\subsection{Estimates on Solutions to the Hyperbolic equation}

Let
  \begin{equation}
   G=K g
   ,
   \llabel{EQ168}
  \end{equation}
where $K$ is a \emph{space/time} linear differential
operator of order
one or two
such that $K$ is tangential when restricted to the boundary.
This means that $K|_{\Gamma_{0} \cup \Gamma_{1}}$ involves derivatives in $x_{1}$, $x_{2}$, and~$t$.
In Lemma~\ref{L05} below, we prove that $Q$ and $Q^2$ satisfy this property.
 
Applying $K$ to the equation \eqref{EQ62} and to the boundary conditions \eqref{EQ63}, we get
  \begin{align}
    &Q^{2}G  
       -  \div_{a}(f \nabla_{a} G) 
       = \tilde{F}
    , \inin \Omega \times (0,T)
    \label{EQ80a}
    \\&
    f\nabla_{a} G \cdot \nu = \tilde{h}  \onon{ \Gamma_{1} \times (0,T)}
    \label{EQ80b}
    \\&
       \partial_{3} G  =\tilde{h}_{0} \onon{ \Gamma_{0} \times (0,T)}
       ,
  \label{EQ80c}
  \end{align}
where
  \begin{align}
      &\tilde{F} = KF
           - [\Delta_{a,f},K]g + [Q^{2},K]g\label{EQ81a}
	 \\&
          \tilde{h}  = Kh +    [f\nabla_{a} ( \cdot )\cdot \nu , K] g
           \label{EQ81b}
	   \\&
           \tilde{h}_{0}   =     [\partial_{3} , K] g
  ,
  \label{EQ81c}
  \end{align}
using the notation
  \begin{equation}
   \Delta_{a,f} = \div_{a} ( f \nabla_{a} \cdot) 
   .
   \llabel{EQ82}
  \end{equation}
In \eqref{EQ81b}, we used the notation
  \begin{equation}
   [f\nabla_{a} ( \cdot )\cdot \nu , K] g
   = f a_{ki} \partial_{k} (K g) \nu_i
       - K (f a_{ki}\partial_{k}g \nu_i)
    .
   \llabel{EQ142}
  \end{equation}
To prove that \eqref{EQ80b} and \eqref{EQ80c} hold, we make essential use of the
fact that $K$ is tangential at the boundary.
Note also that the function $\tilde h$ is defined in the whole $\Omega$---see the comment immediately after Theorem~\ref{T02}.

\cole
\begin{Lemma}
\label{L05}
The operators $Q$ and $Q^2$ are time-tangential at~$\Gamma_0\cup \Gamma_1$.
\end{Lemma}
\colb

The main ingredient in the proof is
  \begin{align}
  \begin{split}
   \vv_i a_{3i}=0
   \onon{\Gamma_0\cup \Gamma_1}
   .
  \end{split}
   \label{EQ138}
  \end{align}

\begin{proof}[Proof of Lemma~\ref{L05}]
Note that $Q=\partial_{t}+\vv_i a_{ki}\partial_{k}=\partial_{t}+\sum_{k=1}^{2}\vv_i a_{ki}\partial_{k}+ \vv_i a_{3i}\partial_{3}$.
The last term vanishes when restricted to $\Gamma_0\cup \Gamma_1$ due to \eqref{EQ138}, and the conclusion for $Q$ is proven.
To obtain the same for $Q^2$, we write it out as
  \begin{align}
  \begin{split}
   Q^2
   &=
   \partial_{tt}
    + \vv_i a_{ki} \partial_{kt}
    + \partial_{t} (\vv_i a_{ki}) \partial_{k}
    + \vv_i a_{ki} \partial_{tk}
    + \vv_i a_{ki}     \partial_{k} (\vv_j a_{mj})\partial_{m}
    + \vv_i a_{ki}     \vv_j a_{mj}\partial_{km}
   .
  \end{split}
   \label{EQ139}
  \end{align}
Denote by $H$ the sum of the terms which are not time-tangential, i.e.,
  \begin{align}
  \begin{split}
   H
   &=
     \vv_i a_{3i} \partial_{3t}
    + \partial_{t} (\vv_i a_{3i}) \partial_{3}
    + \vv_i a_{3i} \partial_{t3}
    + \vv_i a_{ki}     \partial_{k} (\vv_j a_{3j}) \partial_{3}
    + \vv_i a_{ki}     \vv_j a_{3j}\partial_{k3}
    + \sum_{m=1}^{2}\vv_i a_{3i}     \vv_j a_{mj}\partial_{3m}
   .
  \end{split}
   \llabel{EQ140}
  \end{align}
When restricted to $\Gamma_0\cup \Gamma_1$, the first three terms and the last one vanish by \eqref{EQ138}, and we obtain,
using \eqref{EQ138} again,
  \begin{align}
  \begin{split}
   H |_{\Gamma_1}
   &=
     \vv_i a_{ki}     \partial_{k} (\vv_j a_{3j}) \partial_{3}
    + \vv_i a_{ki}     \vv_j a_{3j}\partial_{k3}
   =
    \sum_{k=1}^{2} \vv_i a_{ki}     \partial_{k} (\vv_j a_{3j}) \partial_{3}
   = 0,
  \end{split}
  \llabel{EQ141}
  \end{align}
as claimed.
\end{proof}

\cole
\begin{Lemma}
\label{L01}
The solution $G=Kg$ to \eqref{EQ80a}--\eqref{EQ80c} satisfies the energy inequality
  \begin{align}
  \begin{split}
  &
   \Vert \sqrt{f}Q G(t)\Vert^{2}_{L^{2}}
   + \Vert \sqrt{f}\nabla_{a} G \Vert^{2}_{L^2}
  \\&\indeq
  \leq
   \mathcal{P}_0
   + \int_{0}^{t}( \Vert QG\Vert_{L^2}^2
   + \Vert \nabla_{a} G \Vert_{L^{2}}^{2}
   +  C_{\epsilon}\Vert  G_{t} \Vert_{L^{2}}^{2})
    \mathcal{P}
   + \int_{0}^{t} \Vert J^{-2} f^2\tilde{h}_{0} \Vert^{2}_{H^{1/2}(\Gamma_{0})}
   \\& \indeq\indeq
   +\epsilon\int_{0}^{t}\Vert J^{-1}f \tilde{h} \Vert^{2}_{H^{1/2}(\Gamma_{1})}
   + \int_{0}^{t} \Vert \tilde{F}\Vert^{2}_{L^{2}}
   -\int_{0}^{t}\int_{\Gamma_{0} } \frac1{J^2}\tilde{h}_{0}  f^2\partial_{t}G
   +\int_{0}^{t}\int_{ \Gamma_{1}} \frac{1}{J}\tilde{h} f\partial_{t}G
  ,
  \end{split}
  \label{EQ83}
  \end{align}  
for any $\epsilon\in(0,1]$.
\end{Lemma}
\colb

\begin{proof}
First, we multiply the equation \eqref{EQ80a} by $ f QG$ and integrate in space, leading to
\begin{align}
   \underbrace{  
      \int Q^{2}G f QG}_{A}
        = \underbrace{\int \div_{a} ( f \nabla_{a} G ) f QG}_{\text{B}}
     + \underbrace{ \int \tilde{F} f QG}_{D}
    \label{EQ84}
  .
   \end{align}
For the term $A$ on the left-hand side,
we write
  \begin{equation}
  A=\into f Q^2 G Q G
  =\into f \partial_{t}(QG) QG
   +\into f\vv\cdot \nabla_{a}(QG) QG
   ,
   \llabel{EQ131}
  \end{equation}
from where
  \begin{align}
  \begin{split}
   A
     &=
     \frac12
     \into f \partial_{t} ( (QG)^2)
       + \into f\vv \cdot \nabla_a G_t G_t
  + \int f \vv\cdot\nabla_{a}(\vv\cdot\nabla_{a} G)\,G_{t}
    \\&\indeq\indeq
   +   \int f \vv\cdot \nabla_{a} G_{t}  \vv\cdot\nabla_{a} G
    +   \int f \vv\cdot\nabla_{a}(\vv\cdot\nabla_{a} G) \vv\cdot\nabla_{a}G   
   .
  \end{split}
   \label{EQ132}
  \end{align}
Now, we apply the integration by parts formula 
  \begin{align}
  \begin{split}
   \int g \vv\cdot \nabla_{a} f 
   &=
   - \int f \vv\cdot \nabla_{a}  g
   - \into \div(a\vv) f g 
   ,
  \end{split}
  \label{EQ86}
  \end{align}
which uses the fact that the boundary terms vanish due to
  \begin{equation}
   \vv\cdot \nu= b_{3i}v_{i} - b_{33} \psi_{t}=0
   \onon{\Gamma_0\cup\Gamma_1}
   .
   \llabel{EQ87}
  \end{equation}
By \eqref{EQ86}, the second term on the right-hand side of \eqref{EQ132} equals
$-\frac12\into \div(a \vv f) G_t^2$.
Similarly, the sum of the third and fourth terms in \eqref{EQ132} equals
$   -\into \div(f a \vv) \vv \cdot \nabla_{a} G G_t$, while
the last term in \eqref{EQ132} becomes
$-\frac12 \into \div(f a \vv) (v\cdot\nabla_{a}G)^2$.
Therefore, \eqref{EQ132} leads to
  \begin{align}
  \begin{split}
   A
     &=
     \frac{1}{2}\frac{d}{dt} \Vert \sqrt{f} QG \Vert_{L^2}^{2} 
      - \frac12\int f_t  |  QG |^{2}
      -\frac12\into \div(a \vv f) G_t^2
    \\&\indeq\indeq
      -\into \div(f a \vv) \vv \cdot \nabla_{a} G G_t
      -\frac12 \into \div(f a \vv) (v\cdot\nabla_{a}G)^2
    \\&
     = \frac{1}{2}\frac{d}{dt} \Vert \sqrt{f} QG \Vert_{L^2}^{2} +M_{1}
   .
   \end{split}
   \label{EQ130a}
  \end{align}
The expression $M_1$ is of lower order and satisfies \eqref{EQ83} below.

Next, we consider $B=\int \div_{a} ( f \nabla_{a} G ) f \partial_{t} G
                   +\int \div_{a} ( f \nabla_{a} G ) f \vv\cdot \nabla_{a} G = B_1 + B_2$.
For the first term $B_1$, we have
  \begin{align}
  \begin{split}
  B_{1}
   &= \int \div_{a} ( f \nabla_{a} G )f G_t
   =   \int a_{ki} \partial_{k} (f a_{si} \partial_{s} G) \, f G_t
   \\&
   =  -\int  f a_{si} \partial_{s} G \, f a_{ki}  \partial_{k}  G_t
     - \int  f a_{si} \partial_{s} G \, \partial_{k}(f a_{ki})    G_t
     + \int_{\Gamma_{0} \cup \Gamma_{1}} a_{3i} f a_{si} \partial_{s} G \, f G_t N_3
   \\&
    =  -\frac{1}{2}\frac{d}{dt} \Vert f \nabla_{a} G \Vert_{L^2}^{2}
    + \int f f_t |\nabla_{a}G |^{2} 
    + \int f^{2} \nabla_{a_{t}} G \cdot \nabla_{a}G
    \\& \indeq \indeq
    - \int  f \nabla_{a}G \cdot \div(f a^{T})    G_t
    + \int_{\Gamma_{0} \cup \Gamma_{1}} \frac{1}{J} b_{3i} f a_{si} \partial_{s} G \, f G_t N_3
     \\&
     =  -\frac{1}{2}\frac{d}{dt} \Vert f \nabla_{a} G \Vert_{L^2}^{2}
        + M_{2} 
    - \int_{\Gamma_{0}} \frac{1}{J^2}  f \partial_{3} G \, f G_t
    + \int_{\Gamma_{1}} \frac{1}{J} b_{3i} f a_{si} \partial_{s} G \, f G_t N_3
   ,
  \end{split}
  \llabel{EQ89}
  \end{align}
where we denoted by $(\div c)_i=\partial_{j}c_{ij}$ the divergence of a matrix valued function~$c$.
Using \refer{EQ80b} and \refer{EQ80c}, we thus obtain
  \begin{align}
  \begin{split}
   B_1
    &=
    -\frac{1}{2}\frac{d}{dt} \Vert f \nabla_{a} G \Vert_{L^2}^{2}
    + M_{2} 
    - \int_{\Gamma_{0}} \frac{1}{J^2} \tilde h_0 f^2 G_t
    +\int_{\Gamma_{1}} \frac{1}{J} \tilde h f G_t
   .
  \end{split}
   \llabel{EQ134}
  \end{align}
Note that on $\Gamma_{0}$ we have $\nabla_{a} G \cdot \nu= \partial_{3}G= \partial_{3}K g$.
Next, using Green's theorem, we may rewrite $B_{2}$ as
  \begin{align}
  \begin{split}
   B_{2} &=  \int \div_{a} ( f \nabla_{a} G ) f \vv\cdot \nabla_{a} G
      =  \int a_{ki} \partial_{k}( f a_{si} \partial_{s} G) f \vv_l a_{ml} \partial_{m} G   \\
&=  \underbrace{-\int   f a_{si} \partial_{s} G  \vv_l
        a_{ki} \partial_{k}(f a_{ml} \partial_{m} G) }_{{B}_{21}}
       \underbrace{-\int   f a_{si} \partial_{s} G  \partial_{k} (a_{ki} \vv_l)
      f a_{ml} \partial_{m} G}_{{B}_{22}}  \\
  & \indeq\indeq
    + \underbrace{\int_{\Gamma_0\cup\Gamma_1} a_{3i}  f a_{si} \partial_{s} G f \vv_l a_{ml} \partial_{m} G N_3}_{I_{\text{B}\tau}}
     .
    \end{split}
    \label{EQ90}
  \end{align}
The first term is
of high order but enjoys a hidden symmetry. To realize this, we write
  \begin{align}
  \begin{split}
   f a_{ki} a_{si} \partial_{s} G \partial_{k}(f_{ml}\partial_{m} G)
   &= \frac12 a_{ml} \partial_{m}\bigl(
                                 (f a_{si}\partial_{s} G)
                                 (f a_{ki}\partial_{k} G)
                                \bigr)
   \\&\indeq
    + \Bigl(
       f a_{ki} a_{si} \partial_{s}G \partial_{k}(f a_{ml} \partial_{m} G)
         -        f a_{ml} a_{si} \partial_{s}G \partial_{m}(f a_{ki} \partial_{k} G)
      \Bigr)    
   ,
  \end{split}
   \label{EQ133}
  \end{align}
where we used
  \begin{equation}
   f a_{ml} a_{si} \partial_{s}G \partial_{m} (f a_{ki} \partial_{k} G)
   =
   \frac12 a_{ml} \partial_{m}\bigl(
                                 (f a_{si}\partial_{s} G)
                                 (f a_{ki}\partial_{k} G)
                                \bigr)
   .
   \llabel{EQ136}
  \end{equation}
Note that second term in \eqref{EQ133}
is of lower order since the highest order terms
involving $\partial_{km}G$
cancel out.
Thus, using \eqref{EQ133}, we may rewrite $B_{21}$ as
  \begin{align}
  \begin{split}
   B_{21}
   &=
   \underbrace{
   \frac12
   \into \partial_{m}(\tilde v_l a_{ml} )
                                 (f a_{si}\partial_{s} G)
                                 (f a_{ki}\partial_{k} G)
   }_{{B}_{211}}
   \\&\indeq
   \underbrace{
   -
   \into
   \tilde v_l
   \Bigl(
       f a_{ki} a_{si} \partial_{s}G \partial_{k}(f a_{ml}) \partial_{m} G
         -        f a_{ml} a_{si} \partial_{s}G \partial_{m}(f a_{ki}) \partial_{k} G
   \Bigr)
   }_{{B}_{212}}
   ,
   \end{split}
   \llabel{EQ85}
  \end{align}
where we used
$a_{3l}\tilde v_{l}|_{\Gamma_0\cup \Gamma_1}=0$
and expanded the last term in~\eqref{EQ133}.
Hence, the term $B$ may be expressed as 
  \begin{align}
  \begin{split}
  \int \div_{a} ( f \nabla_{a} G )f QG  
  & =  -\frac{1}{2}\frac{d}{dt} \Vert f \nabla_{a} G \Vert_{L^2}^{2} 
  + M_{2}    + B_{211} + B_{212} +B_{22}  \\
  & \indeq \indeq 
  + \underbrace{ \int_{\Gamma_{0} \cup \Gamma_{1} }
     \frac{1}{J} (f \nabla_{a} G \cdot \nu)( f G_t) }_{I_{\text{\BM}}} N_3 
   + I_{\text{B}\tau}
   .
   \end{split}
   \label{EQ92}
  \end{align}
The terms $M_{1}$, $M_{2}$, $B_{211}$, $B_{212}$, and  $B_{22}$ are estimated directly using H\"older and Sobolev inequalities as
  \begin{align}
  \begin{split}
   M_{1} + M_{2}    + B_{211} + B_{212} +B_{22}\lec  ( \Vert \nabla_{a} G \Vert_{L^{2}}^{2}
   +  \Vert  G_{t} \Vert_{L^{2}}^{2})P\left( \Vert v- \eta_{t}\Vert_{H^{3}} , \Vert \partial_{t}a \Vert_{H^{2}},  \Vert a \Vert_{H^{2}}  \right)
  ,
  \end{split}
  \label{EQ93}
  \end{align}
where $P$ is a polynomial in the indicated arguments.

To treat 
$I_{\text{B}\tau}$, defined in \eqref{EQ90},
note that  the operator $\tilde v\cdot \nabla_{a}$ is tangential on the boundary $\Gamma_{0} \cup\Gamma_{1} $; indeed,
  \begin{align}
   \tilde v\cdot \nabla_{a}
   =
   (v-\eta_{t})\cdot \nabla_{a}
   = \sum_{k=1,2} v_{j} a_{kj} \partial_{k} + (v_{j}a_{3j}- \psi_{t}a_{33})\partial_{3}
   = \sum_{k=1,2} v_{j} a_{kj} \partial_{k}   = v_{1}\partial_{1}+v_{2}\partial_{2}
   .
   \llabel{EQ94}
  \end{align}
We first separate $I_{\text{B}\tau}$ into integrals over $\Gamma_0$, denoted by $I_0$
and $\Gamma_1$, denoted by~$I_1$.
For the first integral, we have
  \begin{align}
  \begin{split}
   I_0
   &=
    -   \int_{\Gamma_0} \frac{1}{J^2}  f \partial_{3} G f \vv_l a_{ml} \partial_{m} G
   =
   -   \int_{\Gamma_0} \frac{1}{J^2}  f \partial_{3} G f (v_1\partial_{1}G+v_2\partial_{2}G)
   \\&
   \lec
   \Vert J^{-2} f^2 \partial_{3} G \Vert_{H^{1/2}(\Gamma_{0})}
   \Vert  v \Vert_{H^{2}}
   \Vert \nabla_{\tau} G \Vert_{H^{-1/2}(\Gamma_{0})}
   \\&
   \lec
   \Vert J^{-2} f^2 \tilde h_0 \Vert_{H^{1/2}(\Gamma_{0})}
   \Vert  v \Vert_{H^{2}}
   \Vert G \Vert_{H^{1}}
   .
  \end{split}
   \label{EQ177}
  \end{align}
Next,
  \begin{align}
  \begin{split}
   I_1
   &=
  \int_{\Gamma_1} a_{3i}  f a_{si} \partial_{s} G f \vv_l a_{ml} \partial_{m} G
   =
  \int_{\Gamma_1} \frac{1}{J} f \tilde h (v_1\partial_{1} G+v_2\partial_{2}G)
  \\&
  \lec
  \Vert J^{-1} f \tilde h \Vert_{H^{1/2}(\Gamma_{1})}
  \Vert  v \Vert_{H^{2}}
  \Vert \nabla_{\tau} G \Vert_{H^{-1/2}(\Gamma_{1})}
  \lec
  \Vert J^{-1} f \tilde h \Vert_{H^{1/2}(\Gamma_{1})}
  \Vert  v \Vert_{H^{2}}
  \Vert G \Vert_{H^{1}}
  .
  \end{split}
   \label{EQ181}
  \end{align}
The inequalities \refer{EQ177} and \refer{EQ181} may be summarized as
  \begin{align}
  \begin{split}
   I_{\text{B}\tau}
   \lec
   \Vert J^{-2}f^2 \tilde h_0\Vert^{2}_{H^{1/2}(\Gamma_{0})}
   + \epsilon \Vert J^{-1} f\tilde h\Vert^{2}_{H^{1/2}(\Gamma_{1})}
   + C_{\epsilon}\Vert  v \Vert^{2}_{H^{2}}  \Vert  G \Vert^{2}_{H^{1}}   
   ,
  \end{split}
   \label{EQ96}
  \end{align}
for $\epsilon\in(0,1]$.
After integrating \eqref{EQ84} in time, and utilizing expressions in \eqref{EQ130a} and \eqref{EQ92} for the terms $A$ and $B$ as well as  $ f= q'(R) \geq 1/C$, then applying the estimates in \eqref{EQ93} and \eqref{EQ96}, we deduce~\eqref{EQ83}. Note that the last two terms in \eqref{EQ83} arise from the term $I_{\text{\BM}}$ in~\eqref{EQ92}.
\end{proof}

The next important lemma estimates expressions occurring in the
proof of Theorem~\ref{T03} below.

\cole
\begin{Lemma}
\label{IBM}
Let~$K=Q^{2}$.
\\
(i)
The trace functions 
 $ 
   \tilde{h}_{0}
$ and $\tilde{h}$ defined in \eqref{EQ81c} and  \eqref{EQ81b} with $K=Q^{2}$
satisfy the estimates
  \begin{align}
   \Vert  \tilde{h}_{0}(t) \Vert_{H^{1/2}(\Gamma_{0})}
   \lec
   \sum_{j=0}^{3}
    \Vert \partial_{t}^{3-j} g \Vert_{H^{j}}\mathcal{P}
  \lec
  \mathcal{P}
   \llabel{EQ98}
  \end{align}
and
  \begin{align}
    \begin{split}
   \Vert  \tilde{h} \Vert_{H^{1/2}(\Gamma_{1})}
    &\lec
     \sum_{j=0}^{2} \Vert  \partial^{2-j}_{t} h  \Vert_{H^{1/2+j}(\Gamma_{1})}
    + \mathcal{P}
   \lec
   \Vert \partial_{t}^{4} w \Vert_{H^{0.5}(\Gamma_{1})}
   +
   \mathcal{P}
  ,
  \end{split}
   \label{EQ99}
   \end{align}
for all~$t>0$.\\
(ii)
The function $\tilde{F}$ defined in \eqref{EQ81a} satisfies
  \begin{align}
  \Vert  \tilde{F} \Vert_{L^2}
  \lec
  \mathcal{P}
.
  \label{EQ100}
  \end{align}
(iii)
The integral
$I_{\BM0} = -\int_{0}^{t}\int_{\Gamma_{0} } \frac{1}{J^2}\tilde{h}_{0}  f^2\partial_{t}Q^{2}g$ satisfies the estimate
  \begin{align}
    \begin{split}
    I_{\BM0}
    &\lec
    \epsilon    \Vert \nabla Q^2 g\Vert_{L^2}^2
     + C_{\epsilon} \mathcal{P}_0
     +  C_{\epsilon} \int_{0}^{t} \mathcal{P}\,ds
    ,
  \end{split}
    \llabel{EQ101}
  \end{align}\\
for every $\epsilon\in(0,1]$.
\\
(iv)
The integral
$I_{\BM1} = \int_{0}^{t}\int_{ \Gamma_{1}}   \frac{1}{J} (Q^{2}-\partial_{t}^{2}) h  f  \partial^{3}_{t}g $ satisfies
  \begin{align}
  \begin{split}
   I_{\BM1}
   &\lec
     \epsilon    \Vert \nabla  g_{tt}\Vert_{L^2}^{9/5}
     + C_{\epsilon} \mathcal{P}_0
     + \epsilon \int_{0}^{t} \Vert\partial_{t}^{4}w \Vert_{H^{1}(\Gamma_1)}^2\,ds
     + C_{\epsilon} \int_{0}^{t}\mathcal{P} 
     ,
  \end{split}
   \llabel{EQ158}
  \end{align}
for every $\epsilon\in(0,1]$.
\\
(v)
The integral $I_{\BM2} =\int_{0}^{t}\int_{ \Gamma_{1}}  \frac{1}{J}    [ f\nabla_{a} \cdot \nu ,   Q^{2}] g   f  \partial^{3}_{t}g$ satisfies
  \begin{align}
   I_{\BM2}
   \lec
     \epsilon    \Vert \nabla \partial_{t}^2 g\Vert_{L^2}^{9/5}
     + C_{\epsilon} \mathcal{P}_0
     +  C_{\epsilon} \int_{0}^{t}\mathcal{P} 
   ,
  \llabel{EQ103}
  \end{align}
for every $\epsilon\in(0,1]$.\\
(vi)
The integral  $I_{\BM3} = -\int_{0}^{t}\int_{ \Gamma_{1}}  \frac{1}{J} \tilde{h} (\partial^{3}_{t} - \partial_{t} Q^{2}) g$ satisfies
  \begin{align}
  \begin{split}
    I_{\BM3}
    \lec
    \mathcal{P}_0
      +  \epsilon     \Vert \partial_{t} g\Vert_{H^{2}}^{9/5}
      +  \epsilon\Vert  g\Vert_{H^3}^{2}
     + \epsilon \int_{0}^{t} \Vert\partial_{t}^{4}w \Vert_{H^{1}(\Gamma_1)}^2\,ds
     + C_{\epsilon} \int_{0}^{t}\mathcal{P} 
   ,
  \end{split}
  \llabel{EQ104}
  \end{align}
for every $\epsilon\in(0,1]$.
\end{Lemma}
\colb

\begin{proof}[Proof of Lemma~\ref{IBM}]
(i)
Recall that
  \begin{align}
    \begin{split}
   \tilde{h}_{0}
    =    [ \partial_{3} ,   Q^{2}] g
   ,
  \end{split}
   \llabel{EQ152}
  \end{align}
and a straight-forward estimate shows that
  \begin{align}
   \begin{split}
   \Vert \tilde{h}_{0}(t)\Vert_{H^{1/2}(\Gamma_{0})} \lec \sum_{j=0}^{3}
   P( \Vert \partial_{t}^{j} v \Vert_{H^{3-j}}, \Vert \partial_{t}^{j} w \Vert_{H^{5/2-j}(\Gamma_1)}, \Vert \partial_{t}^{j} q \Vert_{H^{3-j}})
   \Vert \partial_{t}^{3-j} g(t)\Vert_{H^{j}} 
   .
   \end{split}
   \llabel{EQ153}
  \end{align}
The inequality \refer{EQ99} is obtained by using Lemma~\ref{hF}.
 
(ii) The bound \eqref{EQ100} follows immediately from H\"older and Sobolev inequalities and the fact that the commutators in
the definition of $\tilde{F}$ are second order operators in space and time.

(iii)
The integral $I_{\BM0}$ may be written as
$\int_{0}^{t}\int_{\Gamma_0} H \partial_{t} Q^{2} g$, where
$H= -\frac{1}{J^2} f^2\tilde h_0 $; observe that $\tilde h_0$ is naturally defined in $\Omega$ (not only on $\Gamma_1$) by the formula~\eqref{EQ81c}. 
Let $\zeta=\zeta(x_3)$ be a smooth cut-off function, which equals 1 in a neighborhood of $x_3=0$ and
vanishes for $x_3\geq 1/2$.
Then we have
  \begin{align}
  \begin{split}
   &\int_{0}^{t}\int_{\Gamma_0} H \partial_{t} Q^2 g
   =
   -
   \int_{0}^{t}\int \partial_{3}(H \partial_{t}Q^2 g \zeta(x_3))
   \\&\indeq
   =
    -
    \int_{0}^{t}\int H \partial_{t} \partial_{3}Q^2 g \zeta
    -
    \int_{0}^{t}\int \partial_{3}H \partial_{t}Q^2 g \zeta
    -
    \int_{0}^{t}\int H \partial_{t}Q^2 g \zeta'
   \\&\indeq
   =
    -
    \int H \partial_{3} Q^2 g \zeta
    + \int H  \partial_{3} Q^2g \zeta\Big|_{t=0}
    +
     \int_{0}^{t}\int \partial_{t} H \partial_{3}Q^2 g \zeta
   \\&\indeq\indeq
    -
    \int_{0}^{t}\int \partial_{3}H \partial_{t} Q^2 g \zeta
    -
    \int_{0}^{t}\int H \partial_{t} Q^2 g \zeta'
   ,
  \end{split}
   \label{EQ173}
  \end{align}
where we integrated by parts in time in the term
$    -
    \int_{0}^{t}\int H \partial_{t} \partial_{3}Q^2 g \zeta
$.
The first integral may rewritten as
  \begin{align}
  \begin{split}
    -\int H \partial_{3} Q^2 g \zeta
    &\lec
    \epsilon \int ( \partial_{3} Q^2 g \zeta)^2
    + C_{\epsilon} \int H^2 
   = 
    \epsilon \int ( \partial_{3} Q^2 g \zeta)^2
     + C_{\epsilon} H(0)^2
     +  C_{\epsilon} \int_{0}^{t}\int \partial_{t} (H^2)\,ds
   \\&
    \lec
    \epsilon    \Vert \nabla Q^2 g\Vert_{L^2}^2
     + C_{\epsilon} \mathcal{P}_0
     + C_{\epsilon} \int_{0}^{t} \mathcal{P}\,ds
   ,
  \end{split}
   \llabel{EQ179}
  \end{align}
for every $\epsilon\in(0,1]$.
Clearly, the second integral on the far right of \eqref{EQ173} is bounded by the norms of the initial
data, while the sum of the last three integrals in \eqref{EQ173}
is bounded by
$\int_{0}^{t}\mathcal{P}$.
\colb

(iv)
Here the integral is of the form
$\int_{0}^{t}\int_{\Gamma_1} H \partial_{t}^3 g$, where
$H= \frac{1}{J} f (Q^{2}-\partial_{t}^{2})h $;
here the function $h$, defined by \eqref{EQ65a} on the boundary $\Gamma_1$ may be extended to interior using~\eqref{EQ65}.
Let $\zeta=\zeta(x_3)$ be a smooth cut-off function, which equals 1 in a neighborhood of $x_3=1$ and
vanishes for $x_3\leq 1/2$.
Then we have
  \begin{align}
  \begin{split}
   &\int_{0}^{t}\int_{\Gamma_1} H \partial_{t}^3 g
   =
   \int_{0}^{t}\int \partial_{3}(H \partial_{t}^3 g \zeta(x_3))
   =
    \int_{0}^{t}\int H \partial_{t}^3 \partial_{3}g \zeta
    +
    \int_{0}^{t}\int \partial_{3}H \partial_{t}^3 g \zeta
    +
    \int_{0}^{t}\int H \partial_{t}^3 g \zeta'
   \\&\indeq
   =
    \int H \partial_{t}^2 \partial_{3}g \zeta
    - \int H \partial_{t}^2 \partial_{3}g \zeta\Big|_{t=0}
    -
     \int_{0}^{t}\int \partial_{t} H \partial_{t}^2 \partial_{3}g \zeta
    +
    \int_{0}^{t}\int \partial_{3}H \partial_{t}^3 g \zeta
    +
    \int_{0}^{t}\int H \partial_{t}^3 g \zeta'
   .
  \end{split}
   \label{EQ174}
  \end{align}
\colb
For the first term, we write
  \begin{align}
  \begin{split}
    \int H \partial_{t}^2 \partial_{3}g \zeta
    &\lec   
    \epsilon 
    \Vert \partial_{t}^2 \partial_{3}g\Vert_{L^2}^{9/5}
    +
    C_{\epsilon}
     \Vert H\Vert_{L^2}^{9/4}
    \lec
     \epsilon    \Vert \nabla \partial_{t}^2 g\Vert_{L^2}^2
     + C_{\epsilon} \mathcal{P}_0
     + C_{\epsilon}\int_{0}^{t} \mathcal{P}\,ds
  .
  \end{split}
   \label{EQ180}
  \end{align}
The last two terms in \eqref{EQ174} are estimated directly, while we need
to treat third term
$    - \int_{0}^{t}\int \partial_{t} H \partial_{t}^2 \partial_{3}g \zeta$
more carefully.
The terms resulting from the first, third, and fourth terms in \eqref{EQ65} are
estimated directly by $\int_{0}^{t}\mathcal{P}$, using H\"older and Sobolev inequalities,
so it only remains
to bound
$   
    - \int_{0}^{t}\int \partial_{t} (J^{-1}(Q^2-\partial_{t}^2) \psi_{tt}) \partial_{t}^2 \partial_{3}g \zeta
$.
But this term may also be estimated directly since $Q^2-\partial_{t}^{2}$ does not contain three time derivatives as
  \begin{align}
  \begin{split}
   &
   - \int_{0}^{t}\int \partial_{t} H \partial_{t}^2 \partial_{3}g \zeta
   \lec
   \epsilon \int_{0}^{t} \Vert\partial_{t}^{4}\psi \Vert_{H^{1}}^2\,ds
   + 
   C_{\epsilon}
   \int_{0}^{t}\mathcal{P}
   \lec  
   \epsilon \int_{0}^{t} \Vert\partial_{t}^{4}w \Vert_{H^{1}(\Gamma_1)}^2\,ds
   + 
   C_{\epsilon}
   \int_{0}^{t}\mathcal{P}
   .
  \end{split}
   \label{EQ184}
  \end{align}
\\
(v) is treated similarly to (iv) by considering the extension to the interior and bounded the resulting terms directly.
\\
(vi)
Let $\zeta$ be as in (iv).
Then we have
  \begin{align}
  \begin{split}
   I_{\BM3}
   =
   -
   \int_{0}^{t}\int
    \frac{1}{J}\tilde h \partial_{3}(\partial_{t}^{3}-\partial_{t} Q^2) g\zeta
   - 
   \int_{0}^{t}\int
    \partial_{3}\left(\frac{1}{J}\tilde h\right)( \partial_{t}^{3}-\partial_{t} Q^2 g)
   - 
   \int_{0}^{t}\int
    \frac{1}{J}\tilde h( \partial_{t}^{3}-\partial_{t} Q^2 g) \zeta'
   .
  \end{split}
  \llabel{EQ159}
  \end{align}
The third term is estimated directly by
  \begin{equation}
   \int_{0}^{t}\mathcal{P}\,ds
   ,
   \label{EQ160}
  \end{equation}
so we only need to consider the first two terms, which we denote by~$I_1$
and~$I_2$.

For $I_1$, we use the formula \eqref{EQ81b} for $\tilde h$ and write
  \begin{align}
  \begin{split}
   I_1
    &=
   -
   \int_{0}^{t}\int
    \frac{1}{J}
    Q^{2} h
    \partial_{3}(\partial_{t}^{3}-\partial_{t} Q^2) g
    \zeta
   -
   \int_{0}^{t}\int
    \frac{1}{J}
    [f\nabla_{a} ( \cdot )\cdot \nu , Q^2] g
    \partial_{3}\partial_{t}(\partial_{t}^{2}- Q^2) g
    \zeta
    = I_{11} + I_{12}
    .
  \end{split}
   \label{EQ170}
  \end{align}
In the second integral, we integrate by parts in time. The interior terms
are bounded by \eqref{EQ160}, the pointwise in time terms evaluated at $0$ are estimated by
$\mathcal{P}_0$, while the pointwise in time terms
evaluated at $t$ are bounded as in~\eqref{EQ180}, except that the expression
on the right-hand side needs to be replaced by
  \begin{align}
  \begin{split}
        \epsilon     \Vert \partial_{t} g\Vert_{H^{2}}^{9/5}
      +      \epsilon\Vert g\Vert_{H^3}^{2}
     + C_{\epsilon} \mathcal{P}_0
     + C_{\epsilon} \int_{0}^{t} \mathcal{P}\,ds
  ,
  \end{split}
   \llabel{EQ171}
  \end{align}
where $\epsilon\in(0,1]$ is arbitrary.
Thus we only need to treat the first term in~\eqref{EQ170},~$I_{11}$.
Again, observing \eqref{EQ65}, the integral is split into two parts:
the integrals corresponding to all terms except $\psi_{tt}$ in
\eqref{EQ65} are treated the same way as~$I_{12}$. Therefore, we only
need
to bound
  \begin{align}
  \begin{split}
   I'
   &=
   -
   \int_{0}^{t}\int
    \frac{1}{J}
    Q^{2} \psi_{tt}
    \partial_{3}(\partial_{t}^{3}-\partial_{t} Q^2) g
    .
  \end{split}
   \label{EQ172}
  \end{align}
The important fact about $I'$ is that
$\partial_{t}^{3}-\partial_{t}Q^2$ does not contain $\partial_{t}^3$
and is tangential, and thus we can integrate by parts in a tangential
variable.
Thus we may bound \eqref{EQ172} simply by the far right side of~\eqref{EQ184},
and the estimate of $I_1$ is complete.
In the term $I_2$, we do not  integrate by parts and proceed by simply using Cauchy-Schwarz inequality, and thus $I_2$ can also be bounded
by the bound in~\eqref{EQ184}.
\end{proof}

\cole
\begin{Theorem}
\label{T03}
The solution $g$ to the hyperbolic equation \eqref{EQ62} with the boundary conditions  \eqref{EQ63} satisfies 
  \begin{align}
  \begin{split}
  &  \Vert Q^{3}g\Vert^{2}_{L^{2}}
  + \Vert \nabla Q^{2}g \Vert_{L^2}^{2}
  +
  \left\Vert \partial^{4}_{t} w \right\Vert_{L^2(\Gamma_1)}^{2}
  +  \left\Vert \Deltah \partial^{3}_{t} w\right\Vert_{L^2(\Gamma_1)}^{2} 
  + \int_{0}^{t} \left\Vert \bar\partial \partial_{t}^{4} w \right\Vert_{L^{2}(\Gamma_1)}^{2} \, ds
  \\& \indeq
  \leq C_{\epsilon}\mathcal{P}_0
   +   \epsilon  \Vert q_{tt}\Vert_{H^{1}}^{9/5}
   +  \epsilon     \Vert g_{tt}\Vert_{H^{1}}^{9/5}
   +  \epsilon\Vert \partial_{t} g\Vert_{H^2}^{9/5}
   +  \epsilon\Vert  g\Vert_{H^3}^{2}
     +C_{\epsilon}\int_{0}^{t}
  \mathcal{P} \,ds
  ,
  \end{split}
   \label{EQ106}
  \end{align}
for every $\epsilon\in(0,1]$.
\end{Theorem}
\colb

\begin{proof}
Applying Lemma~\ref{L01} with $G=Q^{2}g$, we have 
  \begin{align}
  \begin{split}
   &
   \Vert \sqrt{f} Q^{3}g \Vert^{2}_{L^{2}} + \Vert \sqrt{f}\nabla_{a} Q^{2}g \Vert_{L^2}^{2}
    \\&\indeq
    \leq
    \mathcal{P}_0
  +\int_{0}^{t}
     (
     \Vert  Q^{3}g \Vert_{L^2}^{2}
      +      \Vert  Q^{2}g \Vert_{H^{1}}^{2}
      +  C_{\epsilon}\Vert  \partial_{t}Q^{2} g \Vert_{L^{2}}^{2}
     )
  \mathcal{P}
  + \int_{0}^{t} \Vert J^{-2} f^2 \tilde{h}_{0} \Vert^{2}_{H^{1/2}(\Gamma_{0})}
  \\& \indeq\indeq
  + \epsilon\int_{0}^{t}\left\Vert J^{-1} f \tilde{h} \right\Vert^{2}_{H^{1/2}(\Gamma_{1})} 
  + \int_{0}^{t} \Vert \tilde{F}\Vert^{2}_{L^{2}}
  \underbrace{-\int_{0}^{t}\int_{\Gamma_{0} } \frac{1}{J^2}\tilde{h}_{0}  f^2\partial_{t}Q^{2}g}_{I_{\BM0} }+\underbrace{ \int_{0}^{t}\int_{ \Gamma_{1}}  \frac{1}{J}  \tilde{h}  f\partial_{t}Q^{2}g }_{I_{\BM}}
  .
  \end{split}
  \label{EQ105}
  \end{align}  
The norms of $\tilde{h}_{0}$, $\tilde{h}$, and $\tilde{F}$ are estimated using Lemma~\ref{IBM}~(i)--(ii),
by the quantity
  \begin{align}
  \begin{split}
   &
   C\epsilon
      \int_{0}^{t} \left\Vert \bar\partial \partial_{t}^{4} w \right\Vert_{L^{2}(\Gamma_1)}^{2} \, ds
   +
   \int_{0}^{t}\mathcal{P}\,ds
   .
  \end{split}
   \llabel{EQ169}
  \end{align}
On the other hand,
the term $I_{\BM0}$ was addressed in Lemma~\ref{IBM}~(iii);
see \eqref{EQ137a} below.
The main task is to treat the
term $I_{\text{BM}}$
Using the formula \eqref{EQ81b}, we have
  \begin{align}
  \begin{split}
    \tilde{h}&
      = Q^{2} h  +  [ f\nabla_{a} \cdot \nu ,   Q^{2}] g          
  .
  \end{split}
   \label{EQ167}
  \end{align}
Using \eqref{EQ167}, the
boundary term $I_{\BM}$ is expressed as
  \begin{align}
  \begin{split}
   I_{\BM}
   &=
   \int_{0}^{t}\int_{\Gamma_{1}}  \frac{1}{J} \tilde{h} f  \partial^{3}_{t}g
   -  \int_{0}^{t}\int_{ \Gamma_{1}} \frac{1}{J} \tilde{h} (\partial^{3}_{t} - \partial_{t} Q^{2}) g
   \\&
   =
   \int_{0}^{t}\int_{ \Gamma_{1}}   \frac{1}{J} Q^{2} h  f  \partial^{3}_{t}g
   +   \int_{0}^{t}\int_{ \Gamma_{1}}    \frac{1}{J}  [ f\nabla_{a} \cdot \nu ,   Q^{2}] g  f  \partial^{3}_{t}g
   -   \int_{0}^{t}\int_{ \Gamma_{1}}   \frac{1}{J} \tilde{h} (\partial^{3}_{t} - \partial_{t} Q^{2}) g
  \\&
   =
   \underbrace{ \int_{0}^{t}\int_{ \Gamma_{1}}  \frac{1}{J}  (Q^{2}-\partial_{t}^{2}) h  f  \partial^{3}_{t}g}_{I_{\BM1}} 
   +  \underbrace{ \int_{0}^{t}\int_{ \Gamma_{1}}   \frac{1}{J}  [ f\nabla_{a} \cdot \nu ,   Q^{2}] g   f  \partial^{3}_{t}g }_{I_{\BM2}}
   \\&\indeq
     \underbrace{  -\int_{0}^{t}\int_{ \Gamma_{1}}  \frac{1}{J}\tilde{h} (\partial^{3}_{t} - \partial_{t} Q^{2}) g }_{I_{\BM3}}
   + \underbrace{ \int_{0}^{t}\int_{ \Gamma_{1}}  \frac{1}{J}  \partial_{t}^{2} h  f  \partial^{3}_{t}g  }_{I_{\BM4}}
   .
   \end{split}
   \llabel{EQ137}
  \end{align}
The terms $I_{\BM1}$, $I_{\BM2}$, and $I_{\BM3}$ are estimated using
Lemma~\ref{IBM}~(iv)--(vi), and we obtain
\begin{align}
\begin{split}
&
I_{\BM0}+
I_{\BM1}+ I_{\BM2} + I_{\BM3}
\\&\indeq
 \lec   \mathcal{P}_0
      +  \epsilon \Vert \nabla Q^{2}g\Vert_{L^2}^{2}
      +  \epsilon     \Vert  g_{tt}\Vert_{H^{1}}^{9/5}
      +  \epsilon\Vert g_t\Vert_{H^2}^{9/5}
      +  \epsilon\Vert  g\Vert_{H^3}^{2}
     + \epsilon \int_{0}^{t} \Vert\partial_{t}^{4}w \Vert_{H^{1}(\Gamma_1)}^2\,ds
     + C_{\epsilon} \int_{0}^{t}\mathcal{P} 
    .
    \label{EQ137a}
\end{split}
\end{align}
Thus it remains to treat 
$I_{\BM4}$, for which we exploit the fact that a term of the same but opposite sign appears in the energy estimates of the plate equation. Note that on $\Gamma_{1}$
  \begin{align}
   f\partial^{3}_{t} g= \frac{1}{R} \partial^{3}_{t}q
   - \underbrace{ \frac{1}{R}\partial^{2}_{t}(fR)  \partial_{t}g - 2\frac{1}{R}   
   \partial_{t}(fR)   \partial^{2}_{t}g }_{\lot}
   .
   \label{EQ107}
  \end{align}
We now focus on the estimate of the problematic highest order term in $I_{\BM4}$,
corresponding to the first term on the right-hand side of \refer{EQ107}, 
which equals
  \begin{align}
  \begin{split}
  \int_{\Gamma_{1}}  \frac{1}{J R} \partial^{2}_{t} h \partial^{3}_{t}q& 
  =\underbrace{ \int_{\Gamma_{1}}  \frac{1}{JR}  \partial^{2}_{t}(\partial_{t} b_{3i} v_{i}) \partial_{t}^{3} q}_{I_{\BM4a}}  
    - \underbrace{\int_{\Gamma_{1}} \frac{1}{JR}\partial^{4}_{t}w \partial_{t}^{3} q }_{I_{\BM4b}}
     \\
   & \indeq 
   + \underbrace{\int_{\Gamma_{1}}\frac{1}{JR} \partial_{t}^{3}q\partial^{2}_{t} \sum_{j=1}^{2} \partial_{j} b_{3i} v_{l} a_{jl}  v_{i}   }_{I_{\BM4c}}
   - \underbrace{ \int_{\Gamma_{1}} \frac{1}{JR}     \partial_{t}^{3} q \partial^{2}_{t}  \sum_{j=1}^{2} \partial_{j}w_{t} v_{l} a_{jl}    }_{I_{\BM4d}} 
   ,
  \end{split}
   \label{EQ108}
  \end{align}
where we used the definition \refer{EQ65} of~$h$.
Note that the term $I_{\BM4a}$ is, except for the factors of $1/JR$ the same as
the terms $I_{\text{B}2}$, $I_{\text{B}3}$, and $I_{\text{B}4}$ in~\eqref{EQ43}.
Thus, we may use the analogous estimates in \eqref{EQ45}, \eqref{EQ187}, and \eqref{EQ46}
to bound~$I_{\BM4a}$.
The terms $I_{\BM4c}$ and $I_{\BM4d}$ are addressed using the same approach as $I_{\BM4a}$ (changing the boundary integral into an interior integral using the FTC) so that
  \begin{align}
    \begin{split}
    &
    I_{\BM4a}  + I_{\BM4c} + I_{\BM4d}
    \\&\indeq
   \lec
   \epsilon  \Vert \partial_{t}^{2} q\Vert_{H^{1}}^{9/5}
   + \epsilon \Vert \partial_{t}^{3}  w \Vert^{2}_{H^{2}(\Gamma_1)}   
   + \epsilon  \int_{0}^{t} \Vert \partial_{t}^{4} w \Vert^{2}_{H^{1}(\Gamma_{1})}
   +    \mathcal{P}_0
   + C_{\epsilon}\int_{0}^{t}\mathcal{P}
   .
  \end{split}
   \llabel{EQ109}
  \end{align}
The term $I_{\BM4b}$ is canceled upon adding the equation \refer{EQ105}, integrated in time, to~\eqref{EQ57}.
From the equations
\eqref{EQ57} and 
\eqref{EQ105},
and using that $f$ is bounded from below, we  obtain the estimate
  \begin{align}
  \begin{split}
  &  \Vert Q^{3}g(t)\Vert^{2}_{L^{2}} + \Vert \nabla_{a} Q^{2}g \Vert^{2}
  +
  \left\Vert \partial^{4}_{t} w(t) \right\Vert^{2}
  +  \left\Vert \Deltah \partial^{3}_{t} w(t) \right\Vert^{2} 
  + \int_{0}^{t} \left\Vert \bar\partial \partial_{t}^{4} w \right\Vert^{2} \, ds
  \\& \indeq
  \leq C_{\epsilon}\mathcal{P}_0
  +\int_{0}^{t}
     (
     \Vert  Q^{3}g \Vert_{L^2}^{2}
      +      \Vert  Q^{2}g \Vert_{H^{1}}^{2} +  \Vert  \partial_{t}Q^{2} g \Vert_{L^{2}}^{2}
      + C_{\epsilon}
     )
  \mathcal{P}
  \\& \indeq\indeq
   +   \epsilon  \Vert \partial_{t}^{2} q\Vert_{H^{1}}^{9/5}
      +  \epsilon\Vert \partial^{2}_{t} g\Vert_{H^{1}}^{9/5}
   +  \epsilon\Vert \partial_{t} g\Vert_{H^2}^{9/5}
      +  \epsilon\Vert  g\Vert_{H^3}^{2}
        \\& \indeq\indeq
   + \epsilon \Vert \partial_{t}^{3}  w \Vert^{2}_{H^{2}(\Gamma_1)}   
   + \epsilon  \int_{0}^{t} \Vert \partial_{t}^{4} w \Vert^{2}_{H^{1}(\Gamma_{1})}
   + \epsilon \Vert \partial_{t}^{2}w \Vert^{2}_{H^{2}(\Gamma_1)}
   +   \epsilon    \Vert \nabla Q^2 g\Vert_{L^2}^2
  ,
  \end{split}
   \label{EQ110}
  \end{align}
for every $\epsilon\in(0,1]$.
For the term $  \epsilon \Vert \partial_{t}^{3}  w \Vert^{2}_{H^{2}(\Gamma_1)}   $, we write
  \begin{align}
  \begin{split}
      \Vert \partial_{t}^{3}  w \Vert^{2}_{H^{2}(\Gamma_1)}
      &\lec
     \Vert \partial_{t}^{3}  \Deltah w \Vert^{2}_{L^{2}(\Gamma_1)}
     + \Vert \partial_{t}^{3}  w \Vert^{2}_{L^{2}(\Gamma_1)}
     \lec
     \Vert \partial_{t}^{3}  \Deltah w \Vert^{2}_{L^{2}(\Gamma_1)}
     + \mathcal{P}_0
     + \int_{0}^{t}\mathcal{P}
     ,
  \end{split}
   \label{EQ185}
  \end{align}
using the FTC in $t$ in the last step.
Note that all the terms on the far right side of \eqref{EQ185} can be absorbed if $\epsilon>0$ is sufficiently small.
The terms $\epsilon \Vert \partial_{t}^{2}  w \Vert^{2}_{H^{2}(\Gamma_1)}   $
and
$   \epsilon    \Vert \nabla Q^2 g\Vert_{L^2}^2$
can also be absorbed.
Next, the eighth term may be absorbed
into the fifth term on the left-hand side, modulo
$\int_{0}^{t}\mathcal{P}\,ds$, which controls the $L^2$ part.
Note that the second term on the right-hand side of \eqref{EQ110}
is bounded by $\int_{0}^{t}\mathcal{P}\,ds$, 
Finally, we replace $\Vert \nabla_{a} Q^{2}g \Vert_{L^2}^{2}$
with   $ \Vert \nabla Q^{2}g \Vert_{L^2}^{2}$ since $a$ is close to the identity
by Lemma~\ref{L03}~(v).
We thus obtain
  \begin{align}
  \begin{split}
  &  \Vert Q^{3}g(t)\Vert^{2}_{L^{2}}
  + \Vert \nabla Q^{2}g \Vert_{L^2}^{2}
  +
  \left\Vert \partial^{4}_{t} w(t) \right\Vert_{L^2(\Gamma_1)}^{2}
  +  \left\Vert \Deltah \partial^{3}_{t} w(t) \right\Vert_{L^2(\Gamma_1)}^{2} 
  + \int_{0}^{t} \left\Vert \bar\partial \partial_{t}^{4} w \right\Vert_{L^2(\Gamma_1)}^{2} \, ds
  \\& \indeq
  \leq
   \epsilon  \Vert \partial_{t}^{2} q\Vert_{H^{1}}^{9/5}
     +  \epsilon\Vert \partial^{2}_{t} g\Vert_{H^{1}}^{9/5}
   +  \epsilon\Vert \partial_{t} g\Vert_{H^2}^{9/5}
      +  \epsilon\Vert  g\Vert_{H^3}^{2}
  +  C_{\epsilon}\mathcal{P}_0
  +  C_{\epsilon}\int_{0}^{t} \mathcal{P} \,ds
  .
  \end{split}
   \llabel{EQ48}
  \end{align}
Using \eqref{EQ83} with $G=Qg$, we may also incorporate the term $\Vert Q^{2} g \Vert_{L^{2}}$ into the left hand side.
This adds to the left side of \eqref{EQ110} also the terms
$\Vert Q^2 g\Vert_{L^2}^2+\Vert \nabla Q g\Vert_{L^2}^2$.
Thus 
we obtain~\eqref{EQ106}.
\end{proof}

\subsection{Full regularity of $R$: elliptic estimates}
To recover the full regularity of $g$, we apply the $H^{3}$ elliptic regularity to \eqref{EQ62} with the boundary conditions~\eqref{EQ63}. 
Thus, $g$ satisfies the elliptic problem
  \begin{align}
  \begin{split}
   &\div_{a}( f \nabla_{a} g ) = Q^{2} g -F \inin{\Omega}
   \\&
   f\nabla_{a}   g \cdot \frac{\nu}{J} = \frac{h}{J}  \onon{\Gamma_{1}}
   \\&
   \partial_{3} g = 0 \onon{\Gamma_0}
   .
  \end{split}
   \label{EQ111}
  \end{align}
The uniform elliptic estimate implies the following statement.

\cole
\begin{Lemma}
\label{L04}
The solution $g$ to the above system satisfies the 
estimates\\
\noindent(i)   
$  \Vert g(t) \Vert_{H^{3}} \lec \mathcal{P}_0
    + \Vert Q^{2} g(t)\Vert_{H^{1}}
    + \int_{0}^{t} \mathcal{P}$,\\
\noindent(ii)  
$   \Vert Qg(t) \Vert_{H^{2}} \lec  \mathcal{P}_0+ \Vert Q^{3} g(t)\Vert_{L^{2}}+ \int_{0}^{t} \mathcal{P}
$,
\\
(iii)        
$   \Vert g_t \Vert_{H^{2}}
\lec
   \Vert Q g\Vert_{H^2}
   + (\Vert v\Vert_{H^2}
        + \Vert w_t\Vert_{H^2(\Gamma_1)})
     \Vert w \Vert_{H^3(\Gamma_1)}
     \Vert g \Vert_{H^{3}}
$, and\\
(iv)       
$  \Vert g_{tt} \Vert_{H^{1}}
    \lec
    \Vert Q^{2} g\Vert_{H^{1}}
    +
    \bigl(
     \Vert g\Vert_{H^{3}}
     +
     \Vert g_t\Vert_{H^{2}}
    \bigr)
    \left(
    \mathcal{P}_0
    +
    \int_{0}^{t}    P\left(\sum_{j=0}^{1}\Vert \partial_{t}^{3-j}v\Vert_{H^{j}},
           \sum_{j=0}^{2}\Vert \partial_{t}^{3-j}w\Vert_{H^{j}(\Gamma_1)}
    \right)
    \right)
\,ds
$.
\end{Lemma}
\colb

By replacing $\Vert g_t\Vert_{H^2}$
in (iv) by its upper bound in (iii), we get
  \begin{align}
  \begin{split}
   \Vert g_{tt} \Vert_{H^{1}}
    \lec
    \Vert Q^{2} g\Vert_{H^{1}}
    +
    \bigl(
     \Vert g\Vert_{H^{3}}
     +
     \Vert Qg\Vert_{H^{2}}
    \bigr)
    \left(
    \mathcal{P}_0
    +
    \int_{0}^{t}    P\left(\sum_{j=0}^{1}\Vert \partial_{t}^{3-j}v\Vert_{H^{j}},
           \sum_{j=0}^{2}\Vert \partial_{t}^{3-j}w\Vert_{H^{j}(\Gamma_1)}
    \right)
    \right)
   \,ds
   .
  \end{split}
   \label{EQ156}
  \end{align}
  
\begin{proof}[Proof of Lemma~\ref{L04}]
(i)
Applying an $H^{3}$ elliptic estimate to the system \eqref{EQ111},
we get
  \begin{align}
   \Vert g(t) \Vert_{H^{3}} \lec \Vert Q^{2} g(t)\Vert_{H^{1}}+ \Vert h(t)  \Vert_{H^{3/2}(\Gamma_{1})} + \Vert  F(t) \Vert_{H^{1}} + \Vert g(t)\Vert_{L^{2}(\Gamma_{0} \cup \Gamma_{1})}
   ,
   \llabel{EQ112}
  \end{align}
with the constants independent of the coefficient matrix $a$ and $f=p'$ for a short time interval due to Lemma~\ref{L01}~(iv) and~(vi).
Writing  $h$ in terms of the time integral of $h_{t}$ and using the estimate on $h_{t}$ and $F$ from Lemma~\ref{hF}. For the last term, we may also apply the  trace inequality and rewrite $g$ as a time integral in terms of $g_{t}$, and then estimate using  $\mathcal{P}_{0}+\int_{0}^{t} \mathcal{P}$ to obtain~(i).

(ii) We apply Q to the elliptic system \eqref{EQ111},
to obtain another elliptic system
  \begin{align}
  \begin{split}
   \div_{a}( f \nabla_{a} Qg ) = Q^{3} g - [Q, \div_{a}( f \nabla_{a} )] g +QF \inin{\Omega}\\
   f\nabla_{a} ( Q g) \cdot \nu = Qh - [ Q, f\nabla_{a} ]g \cdot \nu  \onon{\Gamma_{1}} \\
   f\nabla_{a} (  Qg) \cdot \nu = [ Q, f\nabla_{a} ]g \cdot \nu  \onon{\Gamma_{0}}
  \end{split}
   \label{EQ116}
  \end{align}
satisfied by $Qg$, from which we obtain
  \begin{align}
  \begin{split}
   \Vert Qg \Vert_{H^{2}}
   &\lec \Vert Q^{3} g\Vert_{L^{2}}
   +\left(
     \mathcal{P}_0+\int_{0}^{t}\mathcal{P}\,ds
    \right) \Vert  g\Vert_{H^{2}}
    + \left(
     \mathcal{P}_0+\int_{0}^{t}\mathcal{P}\,ds
    \right)\Vert  g_{t}\Vert_{H^{1}}
   \\&\indeq
   + \Vert Qh\Vert_{H^{1/2}(\Gamma_{1})}
   +   \Vert  QF\Vert_{L^{2}}
   +   \Vert Qg\Vert_{L^{2}(\Gamma_{0} \cup \Gamma_{1})}
  .
  \end{split}
  \llabel{EQ117}
  \end{align}
Note that the commutator
in the first equation of \eqref{EQ116} is a  second order operator in space and time with at most one derivative in time, and the constants again are independent of $a$ and its time derivative for short time, due to Lemma~\ref{L01}.
Writing $h$, $F$, $g$, $g_{t}$, and $Qg$ as time integrals of their time derivatives and estimating $h$ and $F$ using Lemma~\ref{hF}, we obtain~(ii).

(iii) We next deduce the regularity of the first time derivative of $g$ directly from definition of $Qg$ since $\partial_{t} g = Qg - \tilde{v} \nabla_{a}g $.
Thus we have
  \begin{align}
  \begin{split}
   \Vert \partial_{t} g \Vert_{H^{2}}
   &\lec
   \Vert Q g\Vert_{H^2}
   + \Vert v-\eta_t\Vert_{H^2}
     \Vert \nabla_{a}g\Vert_{H^{2}}
   \lec
   \Vert Q g\Vert_{H^2}
   + (\Vert v\Vert_{H^2}+ \Vert \eta_t\Vert_{H^2})
     \Vert a \Vert_{H^2}
     \Vert g \Vert_{H^{3}}
  \\&
  \lec
   \Vert Q g\Vert_{H^2}
   + (\Vert v\Vert_{H^2}
        + \Vert w_t\Vert_{H^2(\Gamma_1)})
     \Vert w \Vert_{H^{3}(\Gamma_{1})}
     \Vert g \Vert_{H^{3}}
    .
  \end{split}
   \llabel{EQ95}
  \end{align}

(iv) The regularity of the second time derivative of $g$ can be deduced from regularity of  the $H^{1}$ norm of $Q^{2}g$ using \eqref{EQ139},
which leads to
\begin{align}
  \begin{split}
  \Vert \partial^{2}_{t} g \Vert_{H^{1}}
   &
    \lec \Vert Q^{2} g\Vert_{H^{1}}
     + \Vert \vv_i a_{ki} \partial_{kt} g\Vert_{H^{1}}
     + \Vert \partial_{t} (\vv_i a_{ki}) \partial_{k} g \Vert_{H^{1}}
    \\&\indeq
      + \Vert \vv_i a_{ki} \partial_{tk} g \Vert_{H^{1}}
      + \Vert \vv_i a_{ki}     \partial_{k} (\vv_j a_{mj})\partial_{m} g \Vert_{H^{1}}
      + \Vert \vv_i a_{ki}     \vv_j a_{mj}\partial_{km} g \Vert_{H^{1}}
    .
  \end{split}
   \llabel{EQ97}
  \end{align}
We estimate the second term as
$\Vert \vv_i a_{ki} \partial_{kt} g\Vert_{H^{1}}
\lec \Vert \vv \Vert_{H^2} \Vert a\Vert_{H^2}\Vert\nabla g_t\Vert_{H^{1}}$ and
the third as
$\Vert \partial_{t} (\vv_i a_{ki}) \partial_{k} g \Vert_{H^{1}}
\lec \Vert \partial_{t}(\vv_i a_{ki})\Vert_{H^{1}} \Vert g\Vert_{H^{3}}$,
with other terms bounded analogously.
We thus get (iv),
completing the proof.
\end{proof}

\startnewsection{Divergence estimates}{sec44}

Here we bound  the divergence of the velocity $v$ in terms of the time derivatives of the density~$R$. Denote the variable divergence by 
  \begin{align}
   \div_{a} v = a_{ki} \partial_{k} v_{i}
   .
   \llabel{EQ119}
  \end{align}
From the continuity equation, we have 
  \begin{align}
  \begin{split}
   a_{ki} \partial_{k} v_{i} =-\frac{ \partial_{t} R}{ R} -  \frac{1}{R} v_{i} a_{ji}\partial_{j} R -  \frac{1}{R}\psi_{t} a_{33}\partial_{3} R
   ,
  \end{split}
   \llabel{EQ120}
  \end{align}
from which we obtain the estimate
  \begin{align}
  \begin{split}
   \Vert \div_{a} \partial_{t}^{s} v_{i} \Vert_{H^{2-s}}
   &\lec 
   \Vert \partial_{t}^{s+1} R \Vert_{H^{2-s}}
   +
   \Vert v\Vert_{H^2}
   \Vert a\Vert_{H^{2}}
   \Vert \partial_{t}^{s} \nabla R \Vert_{H^{2-s}}
   + \mathcal{P}_0
   + \int_{0}^{t}\mathcal{P}\,ds
  ,
  \end{split}
   \llabel{EQ121}
  \end{align}
for $s=1,2$, where $P(0)$ is a polynomial in the norms of the initial conditions.
To return to the divergence, we use the identity
  \begin{equation}
   \div f
   = (\delta_{ij}-a_{ij})\partial_{i}f_j + a_{ij}\partial_{i}f_j
   ,
   \llabel{EQ47}
  \end{equation}
leading to, using Lemma~\ref{L03}~(vi),
  \begin{align}
  \begin{split}
   \Vert \div \partial_{t}^{s} v_{i} \Vert_{H^{2-s}}
   &\lec 
   \Vert \partial_{t}^{s+1} R \Vert_{H^{2-s}}
   +
   \Vert v\Vert_{H^2}
   \Vert a\Vert_{H^{2}}
   \Vert \partial_{t}^{s} \nabla R \Vert_{H^{2-s}}
   + \mathcal{P}_0
   + \int_{0}^{t}\mathcal{P}\,ds
  ,
  \end{split}
   \label{EQ188}
  \end{align}
for~$s=1,2$.

\startnewsection{Vorticity estimates}{sec55}
The Eulerian equation satisfied by $\omega = \curl u$ for compressible flow and barotropic pressure is given by
  \begin{align}
   \omega_{t} + u\cdot \nabla \omega - \omega \cdot \nabla u + (\div u)\omega =0
   \inin{\Omega(t)}
   ,
   \label{EQ122}
  \end{align}
where we used
the well-known fact~$\curl (\rho^{-1}\nabla (p(\rho)))=0$.
Changing to the ALE variable
  \begin{align}
  \begin{split}
  \zeta(x,t)= \omega(\eta(x,t),t)
   ,
  \end{split}
  \llabel{EQ123}
  \end{align}
we may express \eqref{EQ122} as
  \begin{align}
   \partial_{t} \zeta_{i} + v_{j} a_{kj} \partial_{k} \zeta_{i} - \psi_{t}a_{33}\partial_{3} \zeta_{i}- \zeta_{l} a_{kl} \partial_{k} v_{i} + (a_{km}\partial_{k}v_{m} )\zeta_{i} =0
  \inin{\Omega}
   .
   \label{EQ124}
  \end{align}
The vorticity may be expressed as the variable curl in terms of the velocity $v$ as
  \begin{align}
   \zeta(x,t) = \epsilon_{ijk} a_{mj} \partial_{m} v_{k} 
   .
   \llabel{EQ125}
  \end{align}
To obtain an estimate for $\zeta$, we first multiply \eqref{EQ124} by $J$ to get 
  \begin{align}
   J \partial_{t} \zeta_{i} + v_{j} b_{kj} \partial_{k} \zeta_{i} - \psi_{t}\partial_{3} \zeta_{i}- \zeta_{l} b_{kl} \partial_{k} v_{i} + (b_{km} \partial_{k} v_{m} )\zeta_{i} =0
   \inin{\Omega}
   .
   \llabel{EQ126}
  \end{align}
We then use the Sobolev extension $\tilde ~\colon H(\Omega) \to H(\Omega_{0})$, where $\Omega_{0} =\mathbb{T}^{2} \times [-1,2]$, with the support in the $x_{3}$ direction contained in $[-1/2,3/2]$ and extend $J$ by $\bar{J}$ so that $1/4 <\bar{J} <2$ for $x_3 < 2$ and $\bar{J}=0$ for $x_3 \geq 2$. We then consider the equation
  \begin{align}
   \bar{J} \partial_{t} \theta_{i} + \tilde{v}_{j} \tilde{b}_{kj} \partial_{k} \theta_{i} - \tilde{\psi}_{t}\partial_{3} \theta_{i}- \theta_{l} \tilde{b}_{kl} 
    \partial_{k} \tilde{v}_{i} + (\tilde{b}_{km} \partial_{k} \tilde{v}_{m} )\theta_{i} =0
   \inin{\Omega}
   .
   \llabel{EQ127}
  \end{align}
Applying $\Lambda= I-\Delta$ to the equation, we obtain an estimate on the $H^2$ norm of $\zeta$. After a straight-forward estimates, we get
  \begin{align}
  \begin{split}
   \Vert \theta(t) \Vert_{H^{2}(\Omega_0)} 
   & \lec
   \Vert \theta(0) \Vert_{H^{2}(\Omega_0)}
     +\int_{0}^{t} \mathcal{P}
   \lec
   \Vert \zeta(0) \Vert_{H^{2}(\Omega_0)}
     +\int_{0}^{t} \mathcal{P}
   \lec
   \mathcal{P}_0
     +\int_{0}^{t} \mathcal{P}
   .
  \end{split}
   \llabel{EQ128}
  \end{align}
By uniqueness, we then also obtain
  \begin{align}
  \begin{split}
   \Vert \zeta(t) \Vert_{H^{2}} 
   & \lec
    \mathcal{P}_0
     +\int_{0}^{t} \mathcal{P}
   .
  \end{split}
   \llabel{EQ189}
  \end{align}
The derivations in this section are explained in detail in~\cite{KT1,KT2}.

Similarly, we also have
  \begin{align}
  \begin{split}
   \Vert \partial_{t}^{s} \theta(t) \Vert_{H^{2-s}} 
   &\lec
    \mathcal{P}_0
    +
    \int_{0}^{t}\mathcal{P}\,ds
   ,
  \end{split}
   \label{EQ129}
  \end{align}
for~$s=1,2$.

\startnewsection{Div-curl}{sec66}
To recover the full regularity of $v$ and its time derivatives, we use the div-curl type estimate
  \begin{align}
  \begin{split}
   \Vert  \partial_{t}^{s} v \Vert_{H^{3-s}}
   &\lec 
   \Vert  \curl \partial_{t}^{s} v\Vert_{H^{2-s}}
   +\Vert   \div \partial_{t}^{s} v \Vert_{H^{2-s}}  
   +  \Vert  \partial_{t}^{s}    v \Vert_{L^{2}}  
   + \sum_{j=1}^{2}\Vert \partial_{j}^{2-s} \partial_{t}^{s} v_3 \Vert_{L^{2}}  
  \\&
  \lec
   \Vert  \curl_a \partial_{t}^{s} v\Vert_{H^{2-s}}
   +\Vert   \div_a \partial_{t}^{s} v \Vert_{H^{2-s}}
   + \epsilon    \Vert  \partial_{t}^{s} v\Vert_{H^{3-s}}
   \\&\indeq
    +
    \epsilon \Vert q_{tt} \Vert_{H^{1}}^{9/10}
    + \mathcal{P}_0
   + \int_{0}^{t}\mathcal{P}\,ds
  ,
  \end{split}
   \llabel{EQ130}
  \end{align}
for $s=0,1,2$,
where we used \eqref{EQ56} in the last step.
Using \eqref{EQ188} and \eqref{EQ129} and absorbing the term
$\epsilon    \Vert  \partial_{t}^{s} v\Vert_{H^{3-s}}$, we get
  \begin{align}
  \begin{split}
   \Vert  \partial_{t}^{s} v \Vert_{H^{3-s}}
   &\lec 
   \Vert \partial_{t}^{s+1} R \Vert_{H^{2-s}}
   +
   \Vert v\Vert_{H^2}
   \Vert a\Vert_{H^{2}}
   \Vert \partial_{t}^{s} \nabla R \Vert_{H^{2-s}}
   +     \epsilon \Vert R_{tt} \Vert_{H^{1}}^{9/10}
   + \mathcal{P}_0
   + \int_{0}^{t}\mathcal{P}\,ds
   ,
  \end{split}
   \llabel{EQ88}
  \end{align}
whence,
using
$   \Vert v\Vert_{H^2}
   \Vert a\Vert_{H^{2}}
   \Vert \partial_{t}^{s} \nabla R \Vert_{H^{2-s}}
   \lec
   \Vert \partial_{t}^{s} \nabla R \Vert_{H^{2-s}}^2
   +
   \Vert v\Vert_{H^2}^2
   \Vert a\Vert_{H^{2}}^2
$ where the last term can be absorbed in
$\mathcal{P}_0+\int_{0}^{t}\mathcal{P}\,ds$,
we get
  \begin{align}
  \begin{split}
   \Vert  \partial_{t}^{s} v \Vert_{H^{3-s}}
   &\lec 
   \Vert \partial_{t}^{s+1} R \Vert_{H^{2-s}}
   +
   \Vert \partial_{t}^{s} \nabla R \Vert_{H^{2-s}}^2
   +     \epsilon \Vert R_{tt} \Vert_{H^{1}}^{9/10}
   + \mathcal{P}_0
   + \int_{0}^{t}\mathcal{P}\,ds
  ,
  \end{split}
   \label{EQ91}
  \end{align}
for $s=0,1,2$,

\startnewsection{Final estimate}{sec77}
We now multiply
the equations in Lemma~\ref{L04}~(i) and~(ii) with a small constant and add to~\eqref{EQ106}. We may use $R$ and $g$ interchangeably, due to the relation~$g=\log R$. We may then use Lemma~\ref{L04}~(iii) and \eqref{EQ156} to replace instances of time derivatives of $g$ on the right hand side (in the polynomial under the integral), thus obtaining the estimate
  \begin{align}
  \begin{split}
  &  \sum_{j=1}^{3} \Vert Q^{3-j} g(t)\Vert^{2}_{H^{j}} 
  +
    \left\Vert \partial^{4}_{t} w(t) \right\Vert^{2} +  \frac{1}{2}   \left\Vert \Deltah \partial^{3}_{t} w(t) \right\Vert_{L^2(\Gamma_1)}^{2} 
  + \int_{0}^{t} \left\Vert \bar\partial \partial_{t}^{4} w \right\Vert_{L^2(\Gamma_1)}^{2} \, ds
  \\&\indeq
  \leq C_{\epsilon}\mathcal{P}_0
   + \epsilon \Vert g_{tt} \Vert^{9/5}_{H^{1}}  
   +  \epsilon\Vert \partial_{t} g\Vert_{H^2}^{9/5}
      +  \epsilon\Vert  g\Vert_{H^3}^{2}
    + \sum_{j=0}^{3}\int_{0}^{t} P(  \Vert \partial_{t}^{j}v\Vert_{H^{3-j}}^{2},\Vert \partial_{t}^{4-j} w \Vert_{H^{j}(\Gamma_{1})}^{2}, \Vert Q^{3-j}g\Vert_{H^{j}} )
  .
  \end{split}
   \llabel{EQ131}
  \end{align}
Now, we use \eqref{EQ156} and Lemma~\ref{L04}~(iii) to bound the second and third terms on the right, and then apply
Young's inequality  thus obtaining 
  \begin{align}
  \begin{split}
  &  \sum_{j=1}^{3} \Vert Q^{3-j} g(t)\Vert^{2}_{H^{j}} 
  +
    \left\Vert \partial^{4}_{t} w(t) \right\Vert^{2} +  \frac{1}{2}   \left\Vert \Deltah \partial^{3}_{t} w(t) \right\Vert_{L^2(\Gamma_1)}^{2} 
  + \int_{0}^{t} \left\Vert \bar\partial \partial_{t}^{4} w \right\Vert_{\Gamma_1}^{2} \, ds
  \\&\indeq
  \leq C_{\epsilon}\mathcal{P}_0
   + \epsilon (     \Vert Q^{2} g\Vert_{H^{1}}^2
    +   \Vert g\Vert_{H^{3}}^2
   +      \Vert Q g\Vert_{H^{2}}^2
              )
   \\&\indeq\indeq
   + C_{\epsilon}\sum_{j=0}^{3}\int_{0}^{t} P(  \Vert \partial_{t}^{j}v\Vert_{H^{3-j}}^{2},\Vert \partial_{t}^{4-j} w \Vert_{H^{j}(\Gamma_{1})}^{2}, \Vert Q^{3-j}g\Vert_{H^{j}} )
  ,
  \end{split}
   \llabel{EQ157}
  \end{align}
for any $\epsilon>0$. Note that the second term on the right can be absorbed into the left-hand side, which after fixing a sufficiently small $\epsilon>0$ leads to
  \begin{align}
  \begin{split}
  &  \sum_{j=1}^{3} \Vert Q^{3-j} g(t)\Vert^{2}_{H^{j}} 
  +
    \left\Vert \partial^{4}_{t} w(t) \right\Vert_{L^2(\Gamma_1)}^{2} +  \frac{1}{2}   \left\Vert \Deltah \partial^{3}_{t} w(t) \right\Vert_{L^2(\Gamma_1)}^{2} 
  + \int_{0}^{t} \left\Vert \bar\partial \partial_{t}^{4} w \right\Vert_{L^2(\Gamma_1)}^{2} \, ds
  \\&\indeq
  \leq \mathcal{P}_0
   + \sum_{j=0}^{3}\int_{0}^{t} P(  \Vert \partial_{t}^{j}v\Vert_{H^{3-j}}^{2},\Vert \partial_{t}^{4-j} w \Vert_{H^{j}(\Gamma_{1})}^{2}, \Vert Q^{3-j}g\Vert_{H^{j}} )
  .
  \end{split}
   \llabel{EQ161}
  \end{align}
Now, we add \eqref{EQ56} and the square of the div-curl estimate \eqref{EQ91} multiplied by a small
constant to obtain
  \begin{align}
  \begin{split}
  &  \sum_{j=0}^{3}
   \bigl(
     \Vert Q^{3-j} g(t)\Vert^{2}_{H^{j}}
     +   \Vert \partial_{t}^{3-j} v(t)\Vert^{2}_{H^{j}}
   \bigr)
     +  \Vert \partial_{t}^{3} R(t) \Vert_{L^{2}}^{2}
   \\&\indeq\indeq
     +  \sum_{j=0}^{3}
          \Vert \partial^{j}_{t}  w\Vert_{H^{5-j}(\Gamma_1)}^{2}  +\Vert \partial^{4}_{t}  w\Vert_{L^{2}(\Gamma_1)}^{2} 
  + \int_{0}^{t} \left\Vert \bar\partial \partial_{t}^{4} w \right\Vert_{L^2(\Gamma_1)}^{2} \, ds
  \\&\indeq
  \leq
    \epsilon \Vert  \partial_{t}^{2} g(t) \Vert_{H^{1}}^{9/5}
    + C_{\epsilon}\mathcal{P}_0
    + C_{\epsilon} \sum_{j=0}^{3}\int_{0}^{t} P(  \Vert \partial_{t}^{j}v\Vert_{H^{3-j}}^{2},\Vert \partial_{t}^{4-j} w \Vert_{H^{j}(\Gamma_{1})}^{2}, \Vert Q^{3-j}g\Vert_{H^{j}} )
  ,
  \end{split}
  \llabel{EQ321}
  \end{align}
for every~$\epsilon>0$.
Now, we use Lemma~\ref{L04}~(iv) to eliminate the first term
on the right-hand side, also fixing an appropriate~$\epsilon>0$.
This leads to
  \begin{align}
  \begin{split}
  &  \sum_{j=0}^{3}
   \bigl(
     \Vert Q^{3-j} g(t)\Vert^{2}_{H^{j}}
     +   \Vert \partial_{t}^{3-j} v(t)\Vert^{2}_{H^{j}}
   \bigr)
     +  \Vert \partial_{t}^{3} R(t) \Vert_{L^{2}}^{2}
   \\&\indeq\indeq
     +  \sum_{j=0}^{3}
          \Vert \partial^{j}_{t}  w\Vert_{H^{5-j}(\Gamma_1)}^{2}  +\Vert \partial^{4}_{t}  w\Vert_{L^{2}(\Gamma_1)}^{2} 
  + \int_{0}^{t} \left\Vert \bar\partial \partial_{t}^{4} w \right\Vert_{L^2(\Gamma_1)}^{2} \, ds
  \\&\indeq
  \leq
   \mathcal{P}_0
   + \sum_{j=0}^{3}\int_{0}^{t} P(  \Vert \partial_{t}^{j}v\Vert_{H^{3-j}}^{2},\Vert \partial_{t}^{4-j} w \Vert_{H^{j}(\Gamma_{1})}^{2}, \Vert Q^{3-j}g\Vert_{H^{j}} )
  .
  \end{split}
   \label{EQ163}
  \end{align}
A standard Gronwall argument then
gives a uniform upper bound for the left-hand side of~\eqref{EQ163}.
Regularity of the other time derivatives of $g$ and hence $R$ then follow
from Lemma~\ref{L04}~(iii) and~(iv).

\colb
\section*{Acknowledgments}
IK was supported in part by the
NSF grant
DMS-2205493.
\v SN was supported by the Czech Science Foundation (GA\v CR) project
22-01591S and by Praemium Academi\ae  \, of  \v S.~Ne\v casov\'a.
The Institute of Mathematics, CAS, is supported by RVO:67985840.
AT is supported by AUS internal grant FRG23-E-S70.
A part of the work was completed when AT was visiting the  Institute of Mathematics at CAS, Prague.
A part of the work was completed when IK and AT were participating in the program ``Mathematical Problems in Fluid Dynamics, part~2'' at SLMath.

\colb

\colb

\ifnum\sketches=1
\newpage
\begin{center}
  \bf   Notes?\rm
\end{center}
\large
\colb
  \begin{align*}
   & v \in H^{3}
   \\&
    q,R\in H^{3}
   \\&
    w \in L_t^\infty H_x^{5}(\Gamma_1)  
   \\&
    w_t \in L_t^\infty H_x^{4}(\Gamma_1)
   \\&
    w_{tt} \in L_t^\infty H_x^{3}(\Gamma_1)   
   \\&
    w_{ttt} \in L_t^\infty H_x^{2}(\Gamma_1)  
   \\&
    w_{tttt} \in L_t^{\infty}L_x^2(\Gamma_1)   \cap L_t^{2}H_x^{1}(\Gamma_1)
   \\&
    \psi \in L_t^\infty H_x^{5.5}
   \\&
    \psi_t \in L_t^\infty H_x^{4.5}
   \\&
    \psi_{tt} \in L_t^\infty H_x^{3.5}
   \\&
    \psi_{ttt} \in L_t^\infty H_x^{2.5}
   \\&
    \psi_{tttt} \in \cap L_t^{2}H_x^{1.5} (\text{watch out: no good $L^\infty$ space here?})
  \end{align*}
\fi

\end{document}